\documentclass[12pt]{article}
\pdfoutput=1
\usepackage[utf8]{inputenc} 
\usepackage[a4paper]{geometry}
\usepackage[activate={true,nocompatibility},final,tracking=true,kerning=true,spacing=true,factor=1100,stretch=10,shrink=10]{microtype}
\microtypecontext{spacing=nonfrench}
\usepackage{graphicx}
\usepackage{tikz}
\usetikzlibrary{trees,shapes,arrows}

\usepackage{paralist}
\usepackage{subfig}
\usepackage{sidecap}
\usepackage{array}

\usepackage{amsthm}
\usepackage{amsmath}
\usepackage{amsfonts}
\usepackage{mathtools}
\usepackage{todonotes}

\usepackage[ruled,vlined]{algorithm2e}

\usepackage{float}
\usepackage{hyperref}
\hypersetup{
	colorlinks=false,
	linkbordercolor=gray,
	citebordercolor=gray,
	urlbordercolor=gray}%
\usepackage{cleveref}

\usepackage{sectsty}
\allsectionsfont{\sffamily\mdseries\scshape\centering\nohang}

\usepackage{abstract}

\usepackage{caption}
\usepackage{sidecap}

\allowdisplaybreaks

\usepackage[switch,modulo]{lineno}
\usepackage{amssymb}
\usepackage[affil-it]{authblk}

\newtheoremstyle{mplain}
  {2\topsep}   
  {\topsep}   
  {\itshape}  
  {0pt}       
  {\scshape} 
  {:}         
  {5pt plus 1pt minus 1pt} 
  {}          

\renewenvironment{proof}[1][\proofname]{{\scshape #1. }}{\qed}

\theoremstyle{mplain}

\newtheorem{thm}{Theorem}[]
\crefname{thm}{Theorem}{Theorems} 
\Crefname{thm}{Theorem}{Theorems}

\newtheorem{lemma}[thm]{Lemma}
\crefname{lemma}{Lemma}{Lemmas} 
\Crefname{lemma}{Lemma}{Lemmas}

\crefname{conj}{Conjecture}{Conjectures} 
\Crefname{conj}{Conjecture}{Conjectures}

\newtheorem{define}[thm]{Definition}
\crefname{define}{definition}{definitions} 
\Crefname{define}{Definition}{Definitions}

\newtheorem{claim}{Claim}
\crefname{claim}{claim}{claims} 
\Crefname{claim}{Claim}{Claims}

\crefname{obs}{observation}{observations} 
\Crefname{obs}{Observation}{Observations}

\newtheorem{prop}[thm]{Proposition}
\crefname{prop}{Proposition}{Propositions} 
\Crefname{prop}{Proposition}{Propositions}

\crefname{cor}{Corollary}{Corollaries} 
\Crefname{cor}{Corollary}{Corollaries}

\crefname{algocf}{Algorithm}{Algorithms} 
\Crefname{algocf}{Algorithm}{Algorithms}

\newcommand{\eps}{\varepsilon}
\newcommand{\E}{\mathbb{E}}

\DeclareMathOperator{\Poi}{Poi}

\DeclareMathOperator{\SAT}{SAT}
\providecommand{\keywords}[1]{\textbf{\textit{Keywords---}} {\small #1}}

\renewcommand{\Pr}{\mathbb{P}}

\title{Biased random $k$-SAT}
	\author{Joel Larsson\footnote{Mathematics Institute, University of Warwick} \\
	\small\texttt{joel.larsson@warwick.ac.uk}
	\and
	Klas Markstr{\"o}m\footnote{Department of Mathematics and Mathematical Statistics, Ume{\aa} Universitet \\ 
	Research supported by a grant from the Swedish Research Council (Vetenskapsr{\aa}det)}\\
	\small\texttt{klas.markström@math.umu.se}}

\begin{document}
\maketitle

\begin{abstract}
The basic random $k$-SAT  problem is: Given a set of $n$ Boolean variables, and $m$ clauses of size $k$ picked uniformly at random from the set of all such clauses on our variables,  is the conjunction of these clauses satisfiable? 

Here we consider a variation of this problem where there is a bias towards variables occurring positive -- i.e. variables occur negated w.p. $0<p< \frac{1}{2}$ and positive otherwise -- and study how the satisfiability threshold depends on $p$.   For $p<\frac{1}{2}$ this model breaks many of the symmetries of the original random $k$-SAT problem, e.g. the distribution of satisfying  assignments in the Boolean cube is no longer uniform. 
 
For any fixed $k$, we find the asymptotics of the threshold as $p$ approaches $0$ or $\frac{1}{2}$. The former confirms earlier predictions based on numerical studies and heuristic methods from statistical physics.

\end{abstract}
\keywords{Random $k$-SAT, random constraint satisfaction problem, phase transition, combinatorial probability}

\newpage

\section{Introduction}
Random $k$-SAT formulas, and their set of satisfying assignments, have become one of the most studied intersection points of combinatorics, computer science and physics.  The basic problem is as follows. Let $x_1,\ldots, x_n$ be a set of Boolean variables.
A $k$-clause is a Boolean formula of the form $z_i\vee \ldots \vee z_k$, where each $z_i$ is either $x_j$ or $\neg x_j$ for some $j$.
A $k$-SAT formula is a Boolean formula of the form $C_1\wedge\ldots \wedge C_m$, where each $C_i$ is a $k$-clause.
The random $k$-SAT problem asks: if we take a $k$-SAT formula $\Phi$ with $m$ clauses on $n$ variable uniformly at random from all such formulas, for which $m,n$ is $\Phi$ satisfiable w.h.p.? For which $m,n$ is $\Phi$ unsatisfiable w.h.p.? Is there a sharp threshold in between?  It turns out that the crucial parameter is the linear density $\alpha:=m/n$ of $\Phi$. 
 
One of the earliest appearances of this problem is in~\cite{CS}, where it was shown that for certain large $m$ a random CNF of the type just described is not satisfiable, and that this is hard to show using the resolution proof system.   Soon thereafter~\cite{chvatal-reed} demonstrated that for $k=2$ there is a threshold from solvable CNFs to unsolvable ones  at the critical density $\alpha_2:=\frac{m}{n}=1$, and they asked if a similar threshold $\alpha_k$ could be identified for all fixed $k$. Around the same time~\cite{KS,MSL,Selman} a series of extensive -- at least for the computational resources of the time -- simulation studies of this problem was begun and evidence for a threshold constant $\alpha_k$ was presented for many small values of $k$.  This in turn sparked the interest of the statistical physics community, seeing similarities between random $k$-SAT and problems in spin-glass theory. A number of heuristic calculations based on such methods~\cite{MZ1,MZ2} provided conjectures both for the values of $\alpha_k$ and the structure of the set of solution to a satisfiable random CNF.  Parts of these calculations could also be made mathematically rigorous~\cite{Tal}, but not the most crucial ones.

The existence of a sharp threshold for random $k$-SAT was proven by Friedgut in~\cite{Fr}, however neither the location of the this threshold, nor that $\alpha_k$ does not asymptotically depend 
on $n$, follows from his very general threshold results. On the other hand an even more detailed description ~\cite{BBCKW} of the threshold for $k=2$ has been obtained and the asymptotic behaviour of $\alpha_k$ as $k\to\infty$ has been found~\cite{C-O}. The existence of such a constant $\alpha_k$ for large enough $k$ has been established~\cite{DSS},  and for $k\approx\log n$ the exact location of the satisfiability threshold has been determined~\cite{FW,LGX}.  It is generally believed that $\alpha_k$ exists for all fixed $k$ but recently the  exact value given by the cavity-method~\cite{MZ1,MZ2}  has been questioned \cite{LM19a} for low $k$.
 
The aim of this paper is to study a variation on the random $k$-SAT problem where, instead of taking each clause uniformly at random, we introduce a bias parameter $p$ which determines the probability of a variable being negated or non-negated in each clause. Variables in a clause independently occur non-negated with probability $p$ and negated with probability $1-p$, independently for each variable. For $p=\frac{1}{2}$ we thus get the usual random $k$-SAT problem, and for smaller $p$ we get a higher proportion of non-negated variables in our clauses. For $k=2$ this $2$-SAT distribution has been studied~\cite{Austrin} in connection with the hardness of approximating the maximum number of clauses which can be satisfied in random CNF, and surprisingly evidence was found for the balanced case not being the hardest one.  The case $k=2$ is also covered by the results in~\cite{2sat-prescribed-degrees}, where the  threshold behaviour of a more general family of $2$-SAT distributions was identified.   For $k=3$ the threshold has been studied numerically, albeit to quite low precision ~\cite{skewedksat}.   Our main focus will be on determining the satisfiability threshold as a function of both $k$ and $p$, and where possible confirm behaviour  conjectured in the older literature. We will also give some results  on the distribution of satisfying assignments in the hypercube. In the unbiased  version this distribution is uniform and this in part responsible for the inefficiency of some probabilistic tools in analysing the model, as noted in e.g. \cite{AcMO}.

We find the exact threshold for any $p$ when $k=2$, and for $p$ very close to $\frac{1}{2}$ when $k\geq\log n$, similarly to what is known for the balanced case $p=\frac{1}{2}$. For fixed $k\geq 3$, we show that the threshold is approximately quadratic in $p$ near $\frac{1}{2}$ and scales like $p^{1-k}$ as $p\to 0$, the latter confirming a prediction based on replica symmetry heuristics~\cite{biasedksatconjecture}. While the proof for the case $p\to 0$ is mostly an adaption of known methods, the proof for the case $p\to \frac{1}{2}$ is novel, and may be of independent interest. The reader interested in the latter proof can skip directly to \cref{section:incrbias}, which is largely self-contained.

\subsection{Biased k-SAT}

Let $p\in (0,1)$ be a real number. The \emph{biased random $k$-SAT problem} is a random SAT problem, where the clauses are picked according to the following distribution: Start with a set of $n$ Boolean variables and pick from it a $k$-set $K$ uniformly at random. Then, independently for each $x \in K$, pick the literal $x$ with probability $p$ and the literal $\neg x$ otherwise.\footnote{As the problem is symmetric under $p\mapsto 1-p$, we will assume throughout that $p\leq \frac{1}{2}$.} The $k$ literals $z_1, \ldots, z_k$ thus chosen form a clause $C:= z_1\vee z_2\vee \ldots \vee z_k$, which we call a $p$-biased clause. Let $\Phi^p_k(m,n)$ be the conjunction of $m$ i.i.d. $p$-biased clauses.

Let $\alpha_{k}(p,n):= \inf\{\frac{m}{n}:\Pr(\phi^p_k(m,n)\textrm{ is satisfiable})\leq \frac{1}{2}\}$, and let $\alpha_k(p):= \lim_n\alpha_k(p,n)$ if the limit exists. In a slight abuse of notation, we will write (for instance) `$\alpha_3(\frac{1}{2})\geq 3$' even if we do not whether $\alpha_3(\frac{1}{2})$ exists, but this should be understood as being short for `$\liminf_n\alpha_3(\frac{1}{2},n)\geq 3$'.
We will study how the satisfiability threshold $\alpha_k(p)$ behaves as a function of $p$. Regions of particular interest are  $p\approx 0$ and $p\approx \frac{1}{2}$. In the latter case, we will often work with the parametrization $p=\frac{1}{2}-b$ for some small positive $b$.

\subsection{Structured coupon collection}
It is worth noting that random SAT problems are examples of structured coupon collector problems \cite{FLM}. Each assignment of true/false to the $n$ Boolean variables is an $n$-word $(\eps_1, \eps_2, \ldots \eps_n)$ in a two-symbol alphabet, and as such can be identified with the vertices of the hyper-cube $\{-1,1\}^n =: \Sigma_n$. (Here, $-1$ denotes false and $1$ denotes true.)
Each $k$-clause $C$ forbids a subset of those vertices, for instance the clause $C := x_1\vee x_2\vee \neg x_3$ forbids any vertex of the form $(-1,-1,1,\eps_4,\eps_5, \ldots \eps_n)$.
The set of vertices forbidden by $C$ forms a $(n$-$k)$-sub-cube of $\Sigma_n$. 

A sub-cube $C' \subseteq \Sigma_n$ can be represented as an $n$-word $(\delta_1, \delta_2, \ldots \delta_n)$ in the alphabet $\{-1,1,\star\}$ in the following way: 
$\delta_i=-1$ iff $x_i=-1$ for all $x\in C'$, $\delta_i=1$ iff $x_i=1$ for all $x\in C'$, $\delta_i = \star$ otherwise. In other words, $\star$ denotes the `free' coordinates of $C'$ (variables that do not occur in the clause), and the dimension of $C'$ is the number of $\star$'s in $C'$. The set of solutions forbidden by the clause $C$ above is the sub-cube $C'=(-1,-1,1,\star,\ldots \star)$. Since there is a bijection between clauses and sub-cubes in this way, we will henceforth identify a clause with its corresponding sub-cube of forbidden vertices.

We will usually think of random $k$-SAT as a process where we add clauses at integer times until the formula is no longer satisfiable. A variation on this is to add clauses at times given by a Poisson process of intensity $1$. At time $m$, the discrete-time model will have $m$ clauses, while the continuous-time model will have $\Poi(m)$ clauses. These two models are closely related\footnote{Cf. the two random graph models $\mathcal{G}_{n,p}$ and $\mathcal{G}_{n,m}$}, and the satisfiability threshold of either of them is within a multiplicative factor $1\pm o(1)$ of the other. We will therefore work with whichever model best suits our needs at any given time, but take care to specify when we switch between them.

\subsection{Structure of paper}
Each section of the paper is largely self-contained.
We begin by adapting known techniques to the biased version of the random $k$-SAT problem in \cref{section:2sat,section:lowerbound,momentmethods}.

In \cref{section:2sat}, we give the exact threshold for the case $k=2$. The remainder of the paper assumes $k\geq 3$.
In \cref{section:lowerbound}, we find a lower bound on $\alpha_k(p)$ for fixed $k$ by analyzing the unit clause propagation algorithm.
In \cref{momentmethods} we use the method of moments to bound $\alpha_k(p)$ . In \cref{section:momentbounds}, we estimate the first two moments of the number of solutions to a biased $k$-SAT formula. This leads to sharp bounds on $\alpha_k(p)$ for $p$ close to $\frac{1}{2}$ and $k\geq K \log n$ for $K$ sufficiently large. We study a variation on the first moment method in \cref{section:upperbound}, and find an slightly sharper upper bound on $\alpha_k(p)$ which for fixed $k$ is within a constant factor of the lower bound from \cref{section:lowerbound}. These results together establish the asymptotic behavior $\alpha_k(p)\sim p^{1-k}$ as $p\to 0$.

Finally, in \cref{section:incrbias},  we investigate the asymptotics of $\alpha_k(p)$ as $p\to \frac{1}{2}$, by studying how the satisfiability of a formula is affected by changing the occurrences of a single variable. We use a novel combination of tools, including Russo's formula and the Kruskal-Katona theorem, to show that the satisfiability threshold is approximately a parabola near $p=\frac{1}{2}$, i.e. $\alpha_k(p)=\alpha_k+\Theta((p-\frac{1}{2})^2)$.

\section{Satisfiability threshold for biased random 2-SAT} 
\label{section:2sat}
In the special case $k=2$, we find the exact value of the threshold. For the classical (unbiased) $2$-SAT problem, the threshold value of $\alpha_2 = 1$ was established by Chvátal \& Reed~\cite{chvatal-reed} by exploiting some of the structure specific to $2$-SAT.

Later Cooper, Frieze \& Sorkin~\cite{2sat-prescribed-degrees} worked with $2$-SAT formulas of prescribed literal degrees and gave a criterion for satisfiability.
Before we state their theorem, we need the following notation: For a $2$-SAT formula $F$, let $d_i^+$ be the number of occurrences of $x_i$ in $F$, and similarly $d_i^-$ be number of occurrences of $\neg x_i$. For any sequence $\mathbf{d} = (d_1^+,d_1^-, \ldots d_n^+,d_n^-)$, let $D_1:=\sum_i (d_i^++d_i^-)$ and $D_2:=\sum_i d_i^+d_i^-$.
\begin{thm}[CFS]
\label{thm:prescribeddegreesequence}
Let $0<\eps <1$ be constant and $n\to \infty$. Let $\mathbf{d}$ be any literal-degree sequence over $n$ variables with $\max_i d_i^\pm  \leq n^{1/11}$ and $D_1$ even, and let $F$ be chosen uniformly at random from all $2$-SAT formulas with degree sequence $\mathbf{d}$.
\begin{enumerate}[(i)]
\item If $2 D_2 < (1-\eps )D_1$, then 
\(\Pr(F \textrm{ is satisfiable})\to 1\)
\item If $2 D_2 > (1+\eps )D_1$, then 
\(\Pr(F \textrm{ is satisfiable})\to 0\)
\end{enumerate}
Both limits are uniform in $n$ (independent of $d$).
\end{thm}

\begin{thm}
\label{thm:2SAT}
The biased $2$-SAT problem with bias $p$ has a sharp satisfiability threshold at $\alpha_2(p) = \frac{1}{4p(1-p)}=\frac{1}{1-4b^2}$.
\end{thm}
\noindent Sumedha, Krishnamurthy \& Sahoo~\cite{biasedksatontrees} first sketched this, and we will give a full proof.

\noindent \begin{proof}
We prove this for the continuous-time case, the discrete-time case follows.
Pick an $\eps>0$.
For a $p$-biased $2$-SAT formula with $\Poi(m)$ clauses, the quantities $D_1$ and $D_2$ are random variables. $\E D_1 = 2m$, while
\begin{align*}
\E D_2 &= \sum_i \E[d_i^+ d_i^-] = \sum_i \E[d_i^+]\E[ d_i^-] 
\\
&= n\cdot \frac{2m }{n}(1-p) \cdot \frac{2m}{n}p = 4p(1-p) \frac{m^2}{n}
\end{align*}
Both $D_1$ and $D_2$ are sharply concentrated around their means. 
Letting $m = {(1-4\eps)\frac{n}{4p(1-p)}}$, we find that $\|\mathbf{d}\|_\infty =O(\log n)$ with high probability, and
\[
2D_2 < (1+\eps)\cdot 2\E D_2 = (1+\eps)(1-4\eps)\cdot 2m < (1-2\eps) \E D_1 <  (1-\eps)D_1,
\]
with high probability, so that (i) from \cref{thm:prescribeddegreesequence} is satisfied w.h.p.
Similarly, letting $m = (1+4\eps)\frac{n}{4p(1-p)}$ gives that (ii) is satisfied w.h.p. The theorem follows.
\end{proof}

\section{Algorithmic lower bound on satisfiability threshold}
\label{section:lowerbound}
In this section we will show how to adapt the work of Chao and Franco~\cite{3satdiffeq-chao} and Achlioptas~\cite{3satdiffeq-achlioptas} to biased $k$-SAT. We will work in discrete time.

We will show that the algorithm `unit clause propagation' (\ref{algorithm:unitclause}) succeeds in finding a satisfying truth assignment with positive probability when $m$ is not too large.
This algorithm is non-backtracking, and straight-forward to analyze. While better lower bounds are known for the non-biased case, this gives a lower bound that (for any fixed $k$) scales correctly with $p$.

\begin{algorithm}[ht]
\SetAlgoRefName{UCP}
\caption{Unit Clause Propagation}
\DontPrintSemicolon
$V_{\textrm{free}} := \{x: \exists C\in \Phi: x\in C \textrm{ or }\neg x \in C\}$ \tcp*{Variables in $\Phi$}
$V_{\textrm{locked}} := \emptyset$ \tcp*{Variables assigned a value}
\While{$\Phi$ contains non-empty clauses}{
  $U := \{C\in \Phi: |C|=1\}$\;
  \eIf{$U$ is non-empty}{
   let $\ell$ be a literal chosen uniformly at random from $\bigsqcup_{C\in U} C$\;
   let $x$ be the variable of $\ell$ 
   }{
   pick $x$ uniformly at random from $V_{\textrm{free}}$\;
   set $\ell :=x $ with probability $1-p$, $\ell :=\neg x $ otherwise\;
  }
  $V_{\textrm{locked}} := V_{\textrm{locked}}\cup \{\ell\}$  \tcp*{Set $\ell$ to true}
  $V_{\textrm{free}} := V_{\textrm{free}}-\{x\}$ \tcp*{$v$ is no longer free}
  \For{$C \in \Phi$}{
  \If{$\ell \in C$}{
   $\Phi := \Phi - C$ \tcp*{$C$ is satisfied and can be removed} 
   }
  \If{$\neg \ell \in C$}{
   $C:= C- \{\neg \ell\}$ \tcp*{$C$ not yet satisfied} \label{algorithm:unitclause}
   }
   }
 }
\For{$x \in V_{\textrm{free}}$}{
  set $\ell :=x $ with probability $1-p$, $\ell :=\neg x $ otherwise \;
  $V_{\textrm{locked}} := V_{\textrm{locked}}\cup \{\ell\}$ \tcp*{Assign random signs to remaining x}
  }
\eIf{$\Phi$ is empty}{
   \KwRet{$V_{\textrm{locked}}$}  \tcp*{$\Phi$ is empty: no clauses left to satisfy}
   }{
   \KwRet{Failed}\tcp*{$\Phi$ contains an empty clause: unsatisfiable}
  }
\end{algorithm}

\noindent \ref{algorithm:unitclause} is given a SAT formula $\Phi$ as input, in the form of a set of subsets of $\{x_1,\neg x_1, x_2, \neg x_2, \ldots ,\neg x_n\}$. It repeatedly tries to satisfy a unit clause (i.e. a clause on a single variable), and if no such clause exists it sets a random variable to a random value. Satisfied clauses and unsatisfied literals are then removed. The algorithm succeeds iff no empty clauses are ever generated. 

 It might seem strange to let $\ell =v $ with probability $1-p$ in the first `else' of the algorithm, rather than with probability $1$ (since this would maximize the expected number of clauses being satisfied). The reason for this choice is to simplify the analysis of the algorithm: it ensures that the dynamics of the number of $i$-clauses, $i>1$, is independent from the number of $1$-clauses.

\begin{thm}
\label{thm:lowerbound}
There exists a $\delta=\delta(k,p)$ such that if $m$ is at most ${\frac{n}{k^2}(2p(1-p))^{1-k}}$, then \ref{algorithm:unitclause} finds a satisfying assignment to $\Phi_k^p(m,n)$ with probability at least $\delta$.

Furthermore, if a satisfying assignment is found, the number of variables set to `FALSE' follows a binomial distribution with parameters $n$ and $p$.
\end{thm}

\noindent The intuition behind this theorem is as follows: An $i$-clause is turned into a $(i\!-\!1)$-clause with probability $2p(1-p)$, and else removed. So of the $k$-clauses, only a fraction of approximately $(2p(1-p))^{k-1}$ survives long enough to be reduced to $1$-clauses. If the rate at which $1$-clauses are being created is strictly less than the rate at which they are dealt with (which is $1$), queuing theory suggests that the queue of $1$-clauses will remain of bounded size. If the queue is of bounded size, the probability of there existing contradicting $1$-clauses at any given time is of order $n^{-1}$, which suggests that the probability of there ever existing such clauses should be of constant order.

\noindent\begin{proof}
To a large extent the proof is essentially \emph{mutatis mutandis} from~\cite{3satdiffeq-achlioptas},   so we will focus on the necessary modifications and how the results follows from a few key lemmas.

For every $i\in \{0,1,\ldots k\}$ and $j\in \{0,1\ldots, n\}$, let $S_i(j)$ be the number of $i$-clauses at time $j$, and $Y(j) := (S_0(j),S_1(j),\ldots S_k(j))$. Our aim is to understand the trajectory of the (time-inhomogeneous) Markov chain $Y$.

Looking at the expected difference $S_{k}(j+1)-S_{k}(j)$, conditioned on $Y(j)$, we see that 
\[
\E[S_{k}(j+1)-S_{k}(j)|Y(j)] = -\frac{k}{n-j}S_{k}(j)
\]
because any $k$-clause contains $k$ of the $n-j$ free variables, and each of them is equally likely to be locked at time $j+1$.

Similarly, for any $0\leq i<k$, 
\[
\E[S_{i}(j+1)-S_{i}(j)|Y(j)] = -\frac{i}{n-j}S_{i}(j)+2p(1-p)\frac{i+1}{n-j}S_{i+1}(j)
\]
where the first term is as before, while the second term counts the expected number of $(i\!+\!1)$-clauses being shrunk to $i$-clauses. (When a variable occurring in an $(i\!+\!1)$-clause is locked, the clause is shrunk to an $i$-clause with probability $2p(1-p)$, and else removed.) Together these difference equations describe the Markov chain $Y$.

The main proof idea is to look at the scaling limit of the \emph{expected} trajectory of this Markov chain (the so-called \emph{liquid model}). However, there are two problems that arise.

First, for $i\geq 2$, $\E S_i$ is of order $n$ and the dynamics of $Y$ is not too sensitive to deviations in $S_i$ of order $o(n)$. But for $i=0$ or $1$, $\E S_i = O(1)$ and the dynamics is sensitive to small deviations. We deal with this problem by looking at the scaling limit of $S_i$ only for $i\geq 2$, and then making sure that $S_2$ is never so large that the influx of $1$-clauses exceeds the rate at which they can be removed. The following lemma (which we state without proof) is a slight modification of lemma 4 in~\cite{3satdiffeq-achlioptas}.
\begin{lemma}
\label{lemma:nosmallclauses}
For any $\delta, \eps >0$, if $t^*\in(0,1)$ is such that $t^*\leq (1-\eps)$ and 
\[
\textrm{w.h.p.} \quad  S_2(tn)<\frac{(1-\delta)(1-t)}{4p(1-p)} \cdot n \textrm{, for all } 0\leq t\leq  t^*
\]
then $S_0(t^*)=S_1(t^*)=0$ with probability at least $\rho = \rho(\eps,\delta)$.
\end{lemma}
\noindent This lemma says that as long the density of the $2$-clauses stays below $1-\delta$ times the satisfiability threshold for $2$-SAT at time $t\leq t^*$, there is a positive probability that there are no $0$- or $1$-clauses at $t=t^*$. So, next we need to bound the expected value of $S_2$, and show that $S_2$ stays close to it.

To see the second problem, let's look at the system of differential equations describing the scaling limit. 
If we let the functions $c_i$ be defined by $c_i(t):=\lim_{n \to \infty}\frac{1}{n}\E S_i(tn)$ for  $i \geq 2$, they will satisfy the following system (which is the scaling limit of the system of difference equations describing $S_i$):
\begin{align*}
\frac{dc_k\!}{dt}&=-\frac{k}{1-t}c_k(t) , &c_k(0)=c
\\
&\,\,\,\vdots &
\\
\frac{dc_i}{dt}&=(2p(1-p))\frac{(i+1)}{1-t}c_{i+1}(t)-\frac{i}{1-t}c_i(t) , &c_i(0)=0
\\
&\,\,\,\vdots &
\\
\frac{dc_2}{dt}&=(2p(1-p))\frac{3}{1-t}c_3(t)-\frac{2}{1-t}c_2(t) , &c_2(0)=0
\end{align*}
The system has the following unique solution:
\[\forall i \in \{2,\ldots ,k\}, \quad c_i(t)=c\cdot \binom{k}{i}(2p(1-p)t)^{k-i}(1-t)^i\]
We want to show that the trajectory of $S_i$ unlikely to deviate much from the trajectory of $c_i n$. There is a theorem by Wormald~\cite{wormalddiffeq} that lets us do precisely that, and in order for it to apply we need the Markov chain to satisfy these two properties:
\begin{enumerate}[(i)]
\item The system of differential equations describing $c=(c_2,\ldots, c_k)$ can be written on the form $c'(t) = f(c(t),t)$ for some Lipschitz continuous function ${f:D\times I \mapsto \mathbb{R}^{k-1}}$ (for some appropriate domain $D \subseteq \mathbb{R}^{k-1}$ and time interval $I$). 
\item Conditioned on the history of $S_2, \ldots , S_k$ up to time $j$ the probability that ${|S_i(j+1)-S_i(j)|>n^{1/5}}$ is at most $o(n^{-3})$.
\end{enumerate}
These properties are largely unaffected by the value of $p$. The first holds for any time interval $I=[0,1-\eps]$ regardless of $p$. For the second one, the increments follow an approximate Skellam distribution, which have exponential tails.

Since the first condition doesn't hold all the way until time $1$, we will only analyze the algorithm on the time interval $[0,1-\eps]$, and then show that the formula remaining at time $1-\eps$ is sparse enough to be satisfied easily. 
Applying Wormald's theorem, we get the following lemma.
\vspace{-1em}
\begin{lemma}
\label{lemma:expectedtrajectory}
For any $\delta, \eps >0$, there exists $\eta=o(1)$ such that with probability at least $1-\eta$, 
\[
 \quad \big|S_i(tn)-c_i(t)n \big|< \delta n \textrm{, for all }  t < (1-\eps).
\]
\end{lemma}
\noindent In order for \cref{lemma:nosmallclauses} to apply, we need that $S_2(tn)$ is sufficiently small for all $t\leq (1-\eps)$, and now \cref{lemma:expectedtrajectory} tells us that $S_2(tn)$ will stay close to $c_2(t)n$ for all such $t$. That is, for any $\delta>0$, $S_2(tn)/n \leq c_2(t)+\delta$ holds w.h.p. So,
\begin{align*}
c_2(t)+\delta &= c\binom{k}{2} (2p(1-p)t)^{k-2}(1-t)^2+\delta
\\
&\leq \left( ck^2 (2p(1-p))^{k-1}(1-t)+\frac{\delta}{\eps} \right)\cdot \frac{1-t}{4p(1-p)} 
\end{align*}
For this to be at most the upper bound for $S_2$ in \cref{lemma:nosmallclauses}, we need the expression within brackets to be bounded away from $1$. We accomplish that by choosing $c \leq (1-\eps')k^{-2}(2p(1-p))^{1-k}$ and $\delta \leq \eps \eps'/2$ for some $\eps'>0$. Thus, for $c$ as above and $t= 1-\eps$, with probability at least $\rho(\eps,\eps'/2)$ we have that there are no clauses smaller than $2$ at time $t$. 

Let $F$ be the formula remaining at time $t=1-\eps$, conditional on there being no clause of size less than $2$, and let $m_t$ be the number of clauses it consists of. To show that $F$ is satisfiable w.h.p., we construct a new formula $\tilde F$ from $F$ by uniformly at random throwing away literals from each clause with more than $2$ literals. Any  assignment satisfying $\tilde F$ also satisfies $F$. The formula $\tilde F$ is a $2$-SAT formula, and it follows the distribution $\Phi_2^p((1-t)n,m_t)$. By \cref{thm:2SAT},  $\tilde F$ is satisfiable w.h.p. if $\frac{4p(1-p) m_t}{(1-t)n}<1$. 

$S_i$ is close to $c_i n$, so it follows that $m_t < n(\delta+\sum_{i=2}^k c_i(t))$ (with high probability). Thus 
\begin{align*}
4p(1-p) \frac{m_t}{\eps n}
&< 
\frac{4p(1-p)}{\eps} \left({\delta + \sum_{i=2}^k c_i(1-\eps)}\right)
\\
&<\frac{4p(1-p)}{\eps} \left(\delta + c \cdot \sum_{i=0}^k  \binom{k}{i}\big(2p(1-p)\big)^{k-i}\eps^i\right)
\\
&= \frac{4p(1-p)}{\eps}\left({\delta +c\cdot \big(2p(1-p)+\eps\big)^{k}}\right)
\\
&= \frac{4p(1-p)}{\eps}\left(\delta +2p(1-p)k^{-2}\cdot \big(\frac{1+\eps}{2p(1-p)}\big)^{k}\right)
\\
&\leq 4p(1-p)\left(\frac{\delta}{\eps}+ \frac{2p(1-p)}{k^2\eps} \left(1+\frac{\eps}{2p(1-p)}\right)^k\right)
\end{align*}
If we let $\eps = 2p(1-p)/k$ and $\delta = \eps/10$, the expression within brackets is at most $\frac{1}{10}+\frac{1}{k}(1+\frac{1}{k})^k<\frac{9}{10}$. The pre-factor is at most $1$, so the entire expression is bounded away from $1$ and thus the condition of \cref{thm:2SAT} is satisfied. So $\tilde F$ is satisfiable w.h.p., and any assignment that satisfies $\tilde F$ also satisfies $F$.

For the `furthermore' part of the theorem, note that the signs assigned to the variable are i.i.d. Bernoulli r.v.'s., so that the number of variables set to `FALSE' is a binomial random variable with $n$ tries and success probability of $p$.
\end{proof}

\section{Method of moments bounds on satisfiability threshold}
\label{momentmethods}
\noindent The earliest proven upper bound ($\alpha_k\leq 2^k \log 2$) was found by applying the first moment method to the number of satisfying assignments.\footnote{In this section, we work in discrete time.} This upper bound has been improved many times, often by using variations on the same method.

\subsection{Vanilla first and second moment methods}
\label{section:momentbounds}
In this section we will estimate the expected number of satisfying assigment in the biased $k$-SAT model, which leads to an upper bound on $\alpha_k(p)$. 
We will also employ the second moment method, but this only gives a non-trivial lower bound when $k$ is logarithmic in $n$.

While the classical random $k$-SAT problem is vertex transitive on the set $\Sigma_n := \{-1,1\}^n$ of solutions, introducing a bias breaks that symmetry. 
However, the biased version is vertex transitive on any fixed weight `layer' of $\Sigma_n$, i.e.  subset where the number of coordinates equal to $1$ is equal to some constant $i$.
The number of such layers is relatively small ($n+1$, compared to $2^n$ vertices in total), so dealing with each layer separately and then applying a union bound only generates small error terms.
First, we will need some notation. 
\begin{enumerate}[(i.)]
\item For any integers $i,r\in [n]$, define $L_i := \{x\in \Sigma_n: \{j:x_j=1\}=i \}$ to be the $i$:th layer of $\Sigma_n$, and take $x,y\in L_i$ such that $\|x - y\|_1 = 2r$ (i.e. the Hamming distance between $x$ and $y$ is $r$).
\item  Let $C \in \mathcal{S}_k$ be a $p$-biased random sub-cube, and let
$Q(i,n) := \Pr(x\in C)$ and $Q^r(i,n) := \Pr({x\in C} \textrm{ and } {y\in C})$. Note that $Q^0=Q$.
\item  Let $Z_{m,i}$ be the (random) number of non-covered vertices in $L_i$ after $m$ clauses (or cubes) have been drawn, and let $Z_m  := \sum_i Z_{m,i}$.
\item  Finally, let $c_{p,x} := \inf \{c: \E Z_{cn,xn} \leq 1 \}$ and $c_p := \sup_x c_{p,x}$.
\end{enumerate}
Applying the first moment method to the random variable $Z_m$, we see that $\alpha_k(p)\leq (1+o(1))c_p$. 
We will begin by estimating $c_p$. 
\begin{claim} For $i\leq k$, 
\(\big(\frac{i-j+1}{k-j+1}\big)^j \leq {\binom{i}{j}}\big/ {\binom{k}{j}} \leq \big(\frac{i}{k}\big)^j\). \label{ineq:binom}
\end{claim}
\noindent \begin{proof}
Expand the binomial coefficients and repeatedly apply the inequality $\frac{a}{b}<\frac{a+c}{b+c}$ (valid for $0<a<b$ and $c>0$).
\end{proof}

\begin{claim}
Define the binary entropy $H(x):=-x \log(x)-(1-x)\log(1-x)$. Then, for $x,y\in (0,1)$, $H(xy)> x H(y)$ and $\log(1-xy)<x\log(1-y)<-xy$.
\label{ineq:entropy}
\label{ineq:log}
\end{claim}
\noindent \begin{proof}
The second part of the lemma follows from the function $z\mapsto \log(1-z)$ being concave and monotonely decreasing. For the first part, note that $H$ is concave and $H(0)=0$, so $H(xy)>xH(y) + (1-x)H(0)= xH(y)$.
\end{proof}
\begin{claim}
For any $p,x\in (0,1)$, $c_{p,x} = -{\frac{1}{n}\log \binom{n}{xn}}/{\log\big(1-Q(xn,n)\big)}$
\end{claim}
\noindent \begin{proof}
Let $i:=xn$. There are $\binom{n}{i}$ vertices in $L_{i}$, each of which fails to be covered by $m$ $p$-biased random sub-cube with probability $(1-Q(i,n))^m$, so $\E Z_{m,i} = \binom{n}{i} (1-Q(i,n))^m$. This equals $1$ precisely when
\[m=-{\log \binom{n}{i}}\Big/{\log\Big(1-Q(i,n)\Big)}.\]
\end{proof}
\begin{claim} For any $p,x\in (0,1)$, let $\eta_p(x):= x(1-p)+(1-x)p$. Then 
\label{estimate:Q}
$Q(xn,n) = (1-o(1))\cdot  \eta_p(x)^k$
\end{claim}
\noindent\begin{proof}
Let $i:=xn$. For a fixed $v \in L_i$, and a $k$-set $I \subset [n]$, there is a unique cube $C$ containing $v$ and whose set of locked variables is precisely $I$. The number of cubes containing $v$ and with precisely $j$ variables locked to $1$ is then $\binom{i}{j}\binom{n-i}{k-j}$. Each such cube appears with probability $p^j(1-p)^{k-j}\binom{n}{k}^{-1}$, so that
\[
Q(i,n) = \sum_{j=0}^k \frac{\binom{i}{j}\binom{n-i}{k-j}}{\binom{n}{k}}p^j(1-p)^{k-j}
\]
We approximate the binomials using \cref{ineq:binom}:
\[
\frac{\binom{i}{j}\binom{n-i}{k-j}}{\binom{n}{k}} = \left(1+O\Big(\frac{k^2}{n}\Big)-O\Big(\frac{k^2}{i}\Big)\right)\cdot\binom{k}{j} \cdot \left(1-\frac{i}{n}\right)^{k-j}\cdot\left(\frac{i}{n}\right)^j,
\]
where the constants implicit in the big-$O$ notation do not depend on $j$.
Summing over all $j\leq k$, we get that
\begin{align*}
Q(i,n) = \left(1+O\Big(\frac{k^2}{n}\Big)-O\Big(\frac{k^2}{i}\Big)\right)\cdot (x(1-p)+p(1-x))^k
\end{align*}
The claim follows.
\end{proof}
\begin{claim}
$Q^r(i,n)\leq \Big(1-\frac{r}{n-k+1}\Big)^{k-1}Q(i,n)$
\label{estimate:Qr}
\end{claim}
\noindent\begin{proof} First, note that $\binom{n}{k}Q^r(i,n) = \binom{n-r}{k}Q(i-\frac{r}{2},n-r)$. But ${Q(i-\frac{r}{2},n-r)} ={ \frac{n}{n-r}Q(i,n)}$, so
\[
{Q^r(i,n)}/{Q(i,n)}=  \frac{n \binom{n-r}{k}}{(n-r)\binom{n}{k}} = \frac{\binom{n-r-1}{k-1}}{\binom{n-1}{k-1}}.
\]
Together with \cref{ineq:binom} this proves the claim.
\end{proof}

Now that we have some good estimates for the probabilities $Q$ and $Q^r$, we can proceed to use the first moment method to get an upper bound on the satisfiability threshold, and the second moment method to get a lower bound.

\begin{prop}[Bounds on first moment threshold]
For any integer $k\geq 3$ and $p\in (0,\frac{1}{2}]$ (with $k$ and/or $p$ possibly depending on $n$), let $x_+:= \min(\frac{1}{2},\frac{p}{(k-1)(1-2p)})$ and $x_-:=\frac{2}{5} x_+$. Then 
\[H(x_+)\eta_p(x_+)^{-k} \leq c_p <  H(x_+)\eta_p(x_-)^{-k}\]
In particular, 
\begin{enumerate}[(1.)]
\item If $p=\frac{1}{2}-\frac{\beta}{2k}$ for some $\beta\leq 1$, then
\[
\log(2)\cdot 2^k  \leq c_p < \log(2)\cdot 2^k e^{\frac{3}{5}\beta}
\]
\item If $k = \omega(1)$ and $p \in (0,1)$ is fixed, then for all sufficiently large $n$
\[
p^{1-k}\frac{\log k}{k} \cdot e^{-1}\leq c_p< p^{1-k}\frac{\log k}{k}  \cdot e^{-\frac{1}{5}}
\]

\end{enumerate}
\end{prop}
\noindent \begin{proof}
For the lower bound, note that $c_{p}\geq c_{p,x}$ by definition. So, in particular, $c_{p}\geq c_{p,x_+}$, and it suffices to estimate $c_{p,x_+}$. 
For the upper bound, we shall find an small interval on which the supremum is attained. Recall that

\[
c_{p,x} = \frac{\log \binom{n}{xn}}{nQ(xn,n)} = (1-o(1))\cdot H(x)\eta(x)^{-k},
\]
where $\eta(x):=p(1-x)+x(1-p)$. 
For the sake of simplicity we will work with $f(x):=H(x)\eta(x)^{-k}$ rather than directly with $c_{p,x}$.

First, note that $f$ is a strictly concave continuous function on $[0,1]$, whence there exists a unique $x_0$ which maximizes $f$. Since $f(0)=f(1)=0$ whereas $f(x)>0$ for any $x\neq 0,1$, we must have $x_0\in (0,1)$. Second, note that $f$ has a continuous derivative on $(0,1)$, so if $x_-<x_+$ are such that $f'(x_-)>0>f'(x_+)$, then $x_-<x_0<x_+$.
Now, for any $x\in (0,1)$,
\begin{align*}
f'(x)\quad &= \quad H'(x)\cdot\eta(x)^{-k}
-H(x)\cdot k\eta(x)^{-k-1}\eta'(x)
\\
&=\quad \eta(x)^{-k-1} \cdot \underbrace{\left(\eta(x)H'(x)-kH(x)(1-2p)\right)}_{=:\Delta(x)}
\end{align*}
The pre-factor $\eta^{-k-1}$ is always positive, so the sign of $f'$ will be the same as the sign of $\Delta$.
Expanding out the definition of $H$, we can rewrite $\Delta$ as
\begin{align}
\Delta(x)=&(p+x(1-2p))(-\log x+\log(1-x)) \nonumber
\\
&+k(1-2p)(x \log x +(1-x)\log(1-x))\nonumber
\\
=& (\log (1-x) - \log x)\cdot \Big(p-(1-2p)(k-1)x\Big) \label{leftlogterm}
\\
&+ \log(1-x)(1-2p)k \label{rightlogterm}
\end{align}
Let $x_+:= \min(\frac 1 2 ,\frac{p}{(k-1)(1-2p)})$ and consider $\Delta(x_+)$. Then the term (\ref{leftlogterm}) equals $0$ (since at least one of the factors in that term equals $0$), while the term (\ref{rightlogterm}) is negative for any $x \in (0,1)$. So $\Delta(x_+)<0$, and thus $f'(x_+)<0$.

Next, we let $x_-:=\frac{2}{5}x_+$ (which is at most $\frac{1}{5}$) and consider $\Delta(x_-)$. Then the term \ref{leftlogterm} equals
\[
(\log(1-x_-)-\log x_-) \cdot\frac{3}{5}p
\]
which is decreasing in $x_-$ and hence at least $\log\big(\frac{1-\frac{1}{5}}{\frac{1}{5}}\big)\cdot\frac{3}{5}p> 0.8p$. The term \ref{rightlogterm}, on the other hand, equals 
\[
\log(1-x_-)\cdot (1-2p)k = \frac{1}{x_-\!\!\!} \log(1-x_-)\cdot \frac{2kp}{5(k-1)}
\]
which is also decreasing in $x_-$ and hence at least ${5\log\left(1-\frac{1}{5}\right)\frac{2k}{5(k-1)}p}$, which for $k\geq 3$ is at least $ -0.7p$. Together this gives us that $\Delta(x_-)>0.8p-0.7p>0$.
It follows that $f'(x_-)>0$, and together with $f'(x_+)<0$ we have that $x_-<x_0<x_+$.

Now that we have an interval $(x_-,x_+)$ which we know contains $x_0$, we can estimate $f(x_0)=H(x_0) \eta(x_0)^{-k}$. To bound the first factor from above, note that $H$ is a strictly increasing function on $[0,\frac{1}{2}]$, and $x_0< x_+ \leq \frac{1}{2}$. Thus $H(x_0)< H(x_+)$. Similarly, to bound the second factor from above, note that $\eta^{-k}$ is decreasing on $[0,1]$, and $x_-<x_0$, whence $\eta(x_0)^{-k}\leq \eta(x_-)^{-k}$. Together these two inequalities gives us that
\[
f(x_0) \leq  H(x_+) \eta(x_-)^{-k}
\]
This proves the upper bound part of the proposition. For the `in particular'-statements, apply the definition of $H$ and the inequality $H(x) \leq x(1-\log x)$.
\end{proof}

Next, we will show that for $k$ growing sufficiently fast and bias sufficiently small, the first moment bound is tight (i.e. $\alpha_k(p)=(1+o(1)) c_p$). 
\begin{thm}
\label{secondmoment-largek}
Assume $k \geq K\cdot (\log_2 n+\omega(1))$ for some $K\geq 1$ and let $\eps>0$ be fixed.
\begin{enumerate}[(i.)]
\item For any $p$, we have that
\[
\alpha_k(p)\geq c_{p,x_*} = (1-o(1))\cdot H(x_*)\eta_p(x_*)^{-k}
\]
where $x_*$ is defined as the smaller of the two roots of the following equation: \[2x^2\!-\!2x\!-\!2^{-1/K}\!+\!1=0.\]
In other words, $\Phi^{\frac{1}{2}-b}_k(m,n)$ is satisfiable w.h.p. for ${m<(1-\eps)c_{p,x_*} n}$.

\item If $|p-\frac{1}{2}|\leq \frac{1-o_K(1)}{\log n}$, then $\alpha_k(p)=(1+o(1))c_p$. In other words,
\[
\Pr(\Phi^{p}_k(m,n) \textrm{ is satisfiable}) = \begin{cases}
o(1), &m > (1+\eps)c_p n
\\
1-o(1), &m < (1-\eps)c_p n
\end{cases}
\]
\end{enumerate}
\end{thm}
\noindent These results are similar to those known for the unbiased $k$-SAT problem. Setting $K=1$ (and hence $x_*=1/2$) in our theorem we recover the lower bound $k\geq \log_2 n +\omega(1)$ which  is known for that case~\cite{FW,LGX}.   We also note that the result remains valid for $K$ which depend on $n$.

\begin{proof}
Let $I_K:=\{x: 2x^2+2x+2^{-1/K}-1\geq 0\}$.
We will use the second moment method to prove that for any $x\in I_K$, there exists a solution in $L_{xn}$ with high probability if $m<(1-\eps)c_{p,x}$.
From this the two parts of the theorem follow by noting that (i.) $x_*\in I_K$, and (ii.) for $p$ sufficiently close to $\frac{1}{2}$, the $x$ that maximizes $c_{p,x}$ lies in $I_K$.
\[
\frac{\E[Z_{m,i}^2]}{\E[Z_{m,i}]^2} = \frac{\sum_{u,v\in L_i} \Pr(u,v \textrm{ not covered at time }m) }{(\sum_{u\in L_i} \Pr(u\textrm{ not covered at time }m))^2}
\]
In order for the event $\{u\textrm{ not covered by time }m\}$ to occur, a Poisson process of intensity $Q(i,n)$ must have had no event on the the time interval $[0,m]$. The probability of this happening is $\exp(-Q(i,n)m)$. Thus the denominator equals $\binom{n}{i}^2\exp(-2Q(i,n)m)$.

Similarly, in order for the event $\{u,v\textrm{ not covered at time }m\}$ to occur, no clause covering $u$ or $v$ can have occurred. By inclusion-exclusion, the total intensity of clauses covering at least one of $u$ and $v$ is 
$2Q(i,n)-Q^r(i,n)$, where $r=\frac{1}{2}d(u,v)$, so this event has probability $\exp((-2Q(i,n)+Q^r(i,n))m)$.

How many pairs $u,v$ have Hamming distance $2r$? Starting from $u$, such a $v$ is uniquely determined by which $r$ $1$'s we flip to $0$'s, and which $r$ $0$'s we flip to $1$'s. So for any $u\in L_i$ there are $\binom{i}{r}\binom{n-i}{r}$ vertices $v$ at distance $2r$ from $u$. It follows that 
\[
\frac{\E[Z_{m,i}^2]}{\E[Z_{m,i}]^2} =\sum_{r=0}^{i} \exp(Q^r(i,n)m)\cdot\binom{i}{r}\binom{n-i}{r}\binom{n}{i}^{-1}
\]
Define $g$ to be the function
\[g(z):=\log\left(\binom{i}{zn}\binom{n-i}{zn}\binom{n}{i}^{-1+(1-\eps)(1-2z)^{k-1}}\right)\]
so that $\frac{\E[Z_{m,i}^2]}{\E[Z_{m,i}]^2}\leq \sum_{r=0}^{i} e^{g(r/n)}$.

\begin{claim}
The function $g$ is concave on $[\frac{\eps}{4} x,x]$.
\end{claim}

\noindent \begin{proof}
The second derivative of $\log \binom{i}{zn}$ is $-(1+o(1))(\frac{n}{z}+\frac{n}{x-z})$, which is at most $-\frac{n}{x}$. Similarly, the second derivative of $\log \binom{n-i}{zn}$ is at most $-\frac{n}{1-x}$.

Using that $k-1$ is at least $\frac{\log n}{-\log(1-2x(1-x))}$, we see that the second derivative of $ (-1+(1-\eps)(1-2z)^{k-1})\log \binom{n}{i}$ is at most
\begin{align*}
&\log\binom{n}{i}\cdot 4(1-\eps)(k-1)(k-2)\exp\left(-\frac{\log(1-\frac{\eps}{4} x)\log n}{\log(1-2x(1-x))}\right)
\\
<
&H(x) n\cdot  4k^2 \exp\left( -\frac{\eps \log n}{8(1-x)} \right)
\leq
4H(x) k^2 n^{1-\frac{\eps}{4}}
\end{align*}
Hence $\frac{d^2}{dz^2}\log g(z) \leq n\left(-\frac{1}{x}-\frac{1}{1-x}+4H(x) k^2n^{-\frac{\eps}{4}}\right)\leq n(-1+o(1))$.
\end{proof}

\noindent We will estimate $g$ for different ranges of $z$:

\begin{enumerate}[\text{Case} I:]
\item $z<\frac{\eps}{4}  x$ 

The factor $-1+(1-\eps)(1-\frac{r}{n})^{k-1}$ is at most $-\eps$, and ${\binom{n}{i}^{-\eps}\leq \binom{n}{\eps i/2}^{-2}}$. On the other hand, $\binom{i}{zn}\binom{n-i}{zn} \leq \binom{n}{2zn}$, which (by assumption) is at most $\binom{n}{\eps i/2}$. Together this gives us that $e^{g(z)} \leq  \binom{n}{\eps i/2}^{-1}$.
It follows that 
\[
\sum_{r=0}^{\frac{\eps}{4}i}e^{g(r/n)} \leq \frac{\eps i}{4} \cdot \binom{n}{\eps i/2}^{-1}\!\!=o(1)
\]

\item $z\geq z_0:= x(1-x)(1-\delta)$

For convenience, let $y:=2x(1-x)$. Using that $k\geq \frac{\log n +\omega(1)}{-\log(1-y)}$, it follows that
\begin{align*}
&\left(1-2z\right)^{k-1} 
\\
\leq &\exp\left(\frac{\log(1-y(1-\delta))}{-\log(1-y)}\cdot (\log n +\omega(1))\right)
\\
\leq &\exp\left(\Big(-1+\frac{\delta y}{-(1-y)\log(1-y)}\Big)\cdot \log n -\omega(1)\right)
\\
\leq &n^{-1+ \delta/\log 2}\cdot e^{-\omega(1)}
\\
=&o(n^{-1+ 2\delta})
\end{align*}
where the last inequality comes from noting that  $y \mapsto \frac{y}{-(1-y)\log(1-y)}$ is an increasing function on $[0,\frac{1}{2}]$, and thus at most $1/\log 2$.
So $\binom{n}{i}^{(1-2z)^{k-1}}$ is at most $ \exp\big(H(x)\cdot o(n^{2\delta})\big)$, which is $1+o(1)$ if we set $\delta:=\frac{1}{\log n}$. It follows that $e^{g(z)} < (1+o(1))\cdot \binom{i}{zn}\binom{n-i}{zn}\binom{n}{i}^{-1}$ for $z$ in this interval, and summing over $r$ we get
\begin{align*}
\sum_{r=z_0 i}^{i} e^{g(r/n)}  &\leq (1+o(1))\cdot \sum_{r=z_0 i}^i \frac{\binom{i}{r}\binom{n-i}{r}}{\binom{n}{i}}
\\
&\leq (1+o(1))\cdot \underbrace{\sum_{r=0}^i \frac{\binom{i}{i-r}\binom{n-i}{r}}{\binom{n}{i}}}_{=1}
\end{align*}

\item $\frac{\eps}{4}  x \leq  z< z_0$

Since $g$ is concave on this interval, its graph lies beneath any of its tangent lines. In particular, this is true for the tangent line at $z=z_0$. In other words, for any $z$ we have that ${g(z)\leq g(z_0)+(z-z_0)g'(z_0)}$. Since $z-z_0<0$, we need to lower bound $g'(z_0)$ and upper bound $g(z_0)$.
\begin{align}
g'(z_0)\, = \,&n\log\left(\frac{(x-z_0)(1-x-z_0)}{z_0^2}\right)\label{case3:1}
\\
&-2(k-1)\log\binom{n}{i}(1-2z_0)^{k-2}\label{case3:2}
\end{align}
Since $1-2z_0> 1-y \geq \frac{1}{2}$, (\ref{case3:2}) is at most 
\[
2(k-1)\log\binom{n}{i}(1-2z_0)^{k-2}\leq (k-1)H(x)n(1-2z_0)^{k-1} 
\]
and $(1-2z_0)^{k-1}=o(n^{-1})$ by the previous case. The fraction in the right-hand side of (\ref{case3:1}) is at least $\frac{1}{1-\delta}$. Thus
\[
g'(z_0) \geq n\log\left(\frac{1}{1-\delta}\right) - o(k) \geq (1-o(1))\frac{n}{\log n}
\]
To bound $e^{g(z_0)}$ from above, recall from the previous case that it is at most $(1+o(1))\binom{i}{z_0n}\binom{n-i}{z_0n}\binom{n}{i}$. The expression $\binom{i}{zn}\binom{n-i}{zn}\binom{n}{i}$ is maximized for $z=x(1-x)$, where it is $O(n^{-\frac{1}{2}})$. Hence 
\[
g(z) \leq O(1)-\frac{1}{2}\log{n}+(z-z_0)\frac{n}{\log n},
\]
and thus
\begin{align*}
\sum_{r=\frac{\eps}{4} i}^{z_0 i} e^{g(r/n)}
&=
O(n^{-\frac{1}{2}}) \sum_{j=0}^{\infty} e^{-(1-o(1))\frac{j}{n}\cdot \frac{n}{\log n}}
\\
&=
O\left( \frac{n^{-\frac{1}{2}}}{1-e^{-(1-o(1))/\log n}}\right)=O\left(\frac{\log n}{\sqrt{n}}\right)
\end{align*}

\end{enumerate}
Together these three cases give us that $\sum_{r=0}^{z_0 i} e^{g(r/n)} \leq 1+o(1)$, which implies that $\frac{\E[Z_{m,i}^2]}{\E[Z_{m,i}]^2}=1+o(1)$. Chebyshev's inequality then gives us that $Z_{m,i}>0$ with probability $1-o(1)$.
\end{proof}

\subsection{Improved first moment method}
\label{section:upperbound}
When the unit-clause algorithm from \cref{section:lowerbound} succeeds in finding a solution, the number of variables set to `false' is concentrated around $pn$. So we might suspect that $i=pn$ is the dominating term in the expected number of solutions. 
Recall that $Q(i,n,k,p)$ is the probability that a $p$-biased $k$-clause covers an arbitrary vertex in $L_i$, and that for $q(x) := Q(xn,n,k,p)$ we have
\[
q(x) = (1+o(1)) \cdot \big(x(1-p)+(1-x)p\big)^k
\]
Furthermore, recall that $Z_{i,n}$ is the number of uncovered vertices in $L_i$, and that
\[
\frac{1}{n}\log\E Z_{pn,n} = H(p)-c\big(2p(1-p)\big)^k
\]
For small $p$, in order for the right hand side to be negative we need $c$ to be of order $\log \frac{1}{p}\cdot  p^{1-k}$. This only differs from the lower bound by a log-factor, and a small improvement on the vanilla first moment method suffices to correct this: the single-flip method (due to Dubois \& Boufkhad~\cite{singleflip}). We will adapt this method to the biased random $k$-SAT model.
\begin{prop}
\label{prop:upperbound-small-p}
For small enough $p$, $\alpha_k(p)\leq 2p^{1-k}\alpha_k({\frac{1}{2}})$.
\end{prop}

\noindent \begin{proof}
We will apply the first moment method to a subset of solutions, which is guaranteed to be non-empty if the set of solutions is non-empty.

Let $C_1,C_2, \ldots$ be the sub-cubes drawn, and let $K_m := \bigcup_{j=1}^m C_j$ be the set of covered vertices at time $m$.
The hyper-cube can be given a lattice ordering as follows: $u \leq v$ if $u_i \leq v_i$ for every coordinate $i$. This is isomorphic to the lattice ordering on $2^{[n]}$ induced by inclusion. Let $M_{m,i}$ be the number of solutions (uncovered vertices) in $L_i$ that are locally minimal w.r.t. this order.\footnote{We might equally well count the number of locally maximal solutions.} Note that $\sum_i M_{m,i}$ is non-zero if and only if the formula $C_1\wedge \ldots \wedge C_m$ is satisfiable. We will bound $\E M_{m,i}$ from above.

Pick any $u \in L_i$, and let $u^1,\ldots ,u^{n-i}$ be the vertices in $L_{i-1}$ adjacent to $u$. The probability that $u$ is a locally minimal solution can then be written as
\[
\Pr(u \notin K_m)\cdot \Pr(u^1,u^2,\ldots u^{n-i} \in K_m| u \notin K_m)
\]
The event $\{u^j \in K_m\}$ is the union over all $t\in [m]$ of the events $\{u^{j} \in C_t\}$.
For any fixed $t$, the events $\{u^{j} \in C_t\}$ and $\{u^{j'} \in C_t\}$ are mutually exclusive conditional on $\{u \notin K_m\}$, because any cube $C_t$ covering both $u^{j}$ and $u^{j'}$ for some $j\neq j'$ must also cover $u$ (by convexity of $C_t$). We can therefore consider this as a balls-and-bins problem, where the balls are clauses and the bins are vertices $u^j$. Let $X_j := |\{t\in m: u^{j} \in C_t \}|$, i.e. the number of `balls' in `bin' number $j$. Dubhashi \& Ranjan~\cite{negativeassoc-ballsandbins} studied negative dependence of balls-and-bins problems, and in particular showed that the vector of the number of balls in each bin satisfies the negative association property (theorem 13 of that paper). 

So $X=(X_1,X_2,\ldots X_{n-i})$ is a negatively associated vector, and proposition 4 from the same paper gives the following inequality:
\[
\Pr\big(X_j\geq 1,\, \forall j\in[n-i]\big|u \notin K_m\big) \leq \prod_{j=1}^{n-i}\Pr(X_j\geq 1|u \notin K_m)
\]
The left hand side is precisely ${\Pr(u^1,u^2,\ldots u^{n-i} \in K_m| u \notin K_m)}$, while each factor on the right hand side is ${\Pr(u^j\in K_m| u\notin K_m)}$, whence
\begin{align*}
&\Pr(u\notin K_m)\cdot \Pr(u^1,u^2,\ldots u^{n-i}\in K_m| u\notin K_m)
\\
\leq\, &\Pr(u\notin K_m)\cdot \prod_{j=1}^{n-i} \Pr(u^j\in K_m| u\notin K_m)
\\
=\, &\Pr(u\notin K_m)\cdot \Pr(u^1\in K_m| u\notin K_m)^{n-i}
\end{align*}
\begin{claim}
\label{singleflipcondint}
The probability of a $p$-biased random cube covering $u^1$ but not $u$ is $\frac{(1+o(1))p}{\eta_p(x)}\cdot\frac{k}{n} q(x)$
\end{claim}

\noindent \begin{proof}
Assume WLOG that $u^1_1=1$, $u_1=-1$ (and thus $u^1_j=u_j$ for $j>1$). 
For a cube to cover $u^1$ but not $u$, it must have a $1$ in the first position. This happens with probability $p\frac{k}{n}$. The other $n-1$ positions can be seen as a $(k-1)$-co-dimensional sub-cube of $\Sigma_{n-1}$, which must cover the vertex $(u_2,u_3,\ldots,u_n)\in \Sigma_{n-1}$.  This happens with probability \[Q(i,n-1,k-1,p) = \frac{(1+o(1))\cdot Q(i,n,k,p)}{p(1-x)+x(1-p)}\] The claim follows.
\end{proof}

Using \cref{singleflipcondint}, we can calculate the conditional probability of $u^1$ not being covered. 
\begin{align*}
\Pr(u^1 \notin K_m| u \notin K_m)
&= \left(1-\frac{(1+o(1))p }{p(1-x)+x(1-p)}\cdot \frac{k}{n} q(x)\right)^{cn}
\\ 
&= \exp\left( -\frac{(1+o(1))p }{p(1-x)+x(1-p)}\cdot ck q(x) \right)
\end{align*}
from which it follows that $u$ is a locally minimal solution with probability at most
\begin{align*}
&\big(1-\Pr(u^1 \notin K_m| u \notin K_m)\big)^{xn}\big(1-q(x)\big)^{cn}
\\
\leq &\exp\left(xn\log (1-\Pr(u^1 \notin K_m| u \notin K_m))
 +cn  \log (1-q(x))\right)
 \\
 \leq &\exp\left(xn\log\left(\frac{(1+o(1))p }{p(1-x)+x(1-p)}\cdot ck q(x)
 \right)- cn \cdot q(x)\right)
\end{align*}
Recall that $M_{c,x} = M_{c,x} (n,p)$ is the number of locally minimal solutions (in layer $xn$) to a $p$-biased random $k$-SAT instance with $cn$ clauses. We can estimate the expected number of minimal solutions to be at most
\begin{align*}
\frac{1}{n}\log \E M_{c,x} &\leq H(x)+x\log \big(ck q(x)\big) -cq(x)
\\
&\leq x \left(1-\log x +\log \big(ck q(x)\big)-cq(x)/x \right)
\\
&=x \left(1+\log k +\log \big(cq(x)/x\big)-cq(x)/x \right)
\\
&< x \left(1+\log k -cq(x)/2x \right)
\\
&=:U
\end{align*}
Let $c=2p^{1-k}$.
We will use the following bound (valid for sufficiently small $p$ and $n=\omega_p(1)$).
\begin{align*}
\frac{q(x)}{p^{k} }
&=
(1+o(1))\left(\frac{x(1-p)+p(1-x)}{p} \right)^k
\\
&\geq
(1-2pk+o(1))\left(\frac{x}{p}+1\right)^k
>
\frac{1}{2}\left(\frac{x}{p}+1\right)^k
\end{align*}
Then $\frac{c}{2x}(x+p)^k > \frac{1}{2}\cdot \frac{p}{x}(\frac{x}{p}+1)^k$. Let $\varphi_k(t):= {(t+1)^k}/{t}$. This function is minimized for $t={1}/{(k-1)}$, giving the inequality
$\varphi_k(t) \geq {(1+\frac{1}{k-1})^k}/{\frac{1}{k-1}}>(k-1)e$. In particular, for $t=x/p$ we have that
\[
\frac{p}{x} \left(\frac{x}{p}+1\right)^k  \geq (k-1)e,
\]
which in turn gives the following bound on $U$: 
\begin{align*}
U = x \left(1+\log k -cq(x)/2x \right)
\leq x \left(1+\log k -\frac{1}{2}(k-1)e \right)
\end{align*}
The right hand side is less than $-x/2$ for $k\geq 3$.
Hence ${\E M_{c,x}\leq e^{-nx/2}}$, which for $x\geq 1/\sqrt{n}$ is at most $e^{-\sqrt{n}/2}$. 
But for smaller $x$, it suffices to note that $M_{c,x} \leq Z_{c,x}$, and
\[\E Z_{c,x} = \binom{n}{xn}e^{-cq(x)} \leq \exp\left( \sqrt{n} \log n - 2pn \right) =\exp\left(-\Omega(n) \right).\]
Thus, by Markov's inequality, 
\[\Pr(\Phi \textrm{ is unsatisfiable}) \leq \sum_{i=0}^n \E M_{c,i/n} \leq \sqrt{n}e^{-\Omega(n)}+ne^{-\sqrt{n}/2}=o(1).\]
In other words, $\alpha_k(p)\leq c = 2p^{1-k}$ for $p$ small enough.
\end{proof}

\section{Satisfiability threshold for bias near 0}
\label{section:incrbias}
In this section we will prove the main theorem of this paper\footnote{From now on, we will work with the continuous-time version of the biased $k$-SAT problem.} , which describes the shape of the threshold near $p=\frac{1}{2}$. We will work with the parametrization $p=\frac{1}{2}-b$ for some small $b$.
\begin{thm}
\label{thm:parabolicthreshold}
For any $k\geq 3$, there exists a constant ${K_k< 2^{8k}}$ such that  for all sufficiently small $b$,
\[1+2kb^2\leq \frac{\alpha_k(\frac{1}{2}-b)}{\alpha_k(\frac{1}{2})} \leq 1+K_k b^2\]
\end{thm}
\noindent The upper bound on $K_k$ can be improved slightly to $2^{6k-4}$ by optimizing the choice of parameters in the proof of \cref{sparsity:hypergraph} due to Chvátal-Szemeredi~\cite{CS}, but in order to achieve a sub-exponential upper bound \cref{lemma:spinesupperbound} would need to be significantly sharpened.

As an aside, what should we expect the correct quadratic coefficient to be?
While the lower bound $1+2kb^2$ matches our exact result of $1+4b^2+O(b^4)$ for $k=2$, the upper and lower bounds for larger bias suggest a larger quadratic coefficient for $k\geq 3$. The threshold curve for $p$ near $0$ or $1$ scales like $(4p(1-p))^{1-k}$, and \emph{if} it follows that curve for $p$ near $1/2$ too, we get the Taylor expansion $1+4(k-1)b^2+O(b^4)$ near $b=0$. 
One might therefore guess that $4(k-1)$ is the correct coefficient. This also matches the result for $k=2$.

Before we can continue with the proof of \cref{thm:parabolicthreshold}, we will need some background on \emph{spine variables}.
\begin{define}
Let $\Phi$ be a satisfiable formula and $x$ a variable in it. We say that $x$ is a \emph{spine variable} in $\Phi$ if $x$ has the same value in any assignment satisfying $\Phi$. If such an $x$ always has value `TRUE', we say that it is a \emph{positive spine variable} and that it is \emph{locked to TRUE}. (Similarly for negative.)

For $\Phi$ an unsatisfiable formula, we say that $x$ is a spine variable in $\Phi$ if there exists a satisfiable formula $F\subset \Phi$ such that $x$ is a spine variable in $F$.
\end{define}

We will use the following definition and lemma from Chvátal-Szemerédi~\cite{CS}.
\begin{define}
Let $x,y>0$. A $k$-uniform hypergraph with $n$ vertices is {\em
($x$,$y$)-sparse} if every set of $s\leq xn$ vertices contains at
most $ys$ edges.
\end{define}

\begin{lemma}[Chvátal-Szemerédi]
\label{sparsity:hypergraph}
Let $k,t>0$ and $y>1/(k-1)$. Then w.h.p.\ the
random $k$-uniform hypergraph with $n$ vertices and $tn$ edges is
$(x,y)$-sparse, where
\begin{equation*}\label{x-sparse}
x = \left( \frac{1}{2e}\left(\frac{y}{te}\right)^{y}\right)^{\frac{1}{y(k-1)-1}}
\end{equation*}
\end{lemma}
\noindent Using this lemma, Boetcher, Istrate \& Percus proved that the number of spine variables is either $o(n)$ or at least $\delta n$ for some $\delta>0$ (part 1 of theorem 3 in~\cite{spinevars}). But their proof actually works even when replacing $o(n)$ with $0$, leading to the following theorem.
\begin{thm}[Boetcher-Istrate-Percus\footnote{A special case of their Theorem 3, using our notation}]
\label{thm:spinevars}
For $k\geq 3$, there is a constant $\delta=\delta(t,k)>0$ such that the following holds with high probability: for every unsatisfiable $\Phi = \Phi_k^{\frac{1}{2}-b}(t'n,n)$ with $t'<t$, the number of spine variables of $\Phi$ is either $0$ or at least $\delta n$.
\end{thm}

\noindent The following lemma is a consequence of the proofs of \cref{sparsity:hypergraph,thm:spinevars}, and we state it without proof.
\begin{lemma} \label{delta-asymptotics}
For a fixed $k$, $\delta=\delta(t,k)$ is a decreasing function of $t$. Furthermore, for a fixed $t<2^k$, $\delta=\exp(-O(k))$. 
\end{lemma}

\noindent Next, we will need a lemma that gives a correspondence between spine variables and clauses that turn a satisfiable formula into an unsatisfiable one. 
\begin{lemma}
\label{spinevar:sat}
Let $\Phi$ be a satisfiable $k$-CNF with a set $S_+\subseteq [n]$ of positive spine variables and a set $S_-\subseteq [n]$ of negative spine variables. Then $\Phi\wedge C$ is unsatisfiable if and only if $C$ can be written as 
\[
C(x) = \Bigg(\bigvee_{i\in K_-} x_i \Bigg)\vee \Bigg(\bigvee_{i\in K_+} - x_i\Bigg)
\]
for some $K_\pm \subseteq S_\pm$.
\end{lemma}
\noindent In other words, $\Phi\wedge C$ is unsatisfiable iff every variable that occurs in $C$ is a spine variable in $\Phi$, and the sign that it has in $C$ is incompatible with the truth value it is locked to in $\Phi$.

\noindent \begin{proof}
For the `if' part, assume $C$ is of the above form. If $x\in \{\pm 1\}^n$ is a solution to $\Phi$, it has $x_i=- 1$ for $i\in S_-$, so ${\bigvee_{i\in K_-} x_i =-1}$.
Similarly, ${\bigvee_{i\in K_+} -x_i =-1}$.
But then $C(x) = -1$, so $x$ is not a solution to $\Phi \wedge C$.

On the other hand, if $x$ is not a solution to $\Phi$ it is not a solutions to $\Phi \wedge C$. Hence $\Phi\wedge C$ is unsatisfiable, proving the `if' part. 

For the `only if' part, assume instead $C$ is not of that form. Then either
\begin{enumerate}
\item There exists an $i\in S_-$ such that $x_i$ occur in $C$ with negative sign. In that case, any $x$ that satisfies $\Phi$ will also satisfy $C$, because such an $x$ will have $x_i=- 1$, and $C$ will have a term $-x_i$. Hence $\Phi \wedge C$ is satisfied by $x$.
\item There exists an $i\in S_+$ such that $x_i$ occur in $C$ with positive sign. Analogously to the previous case, any $x$ that satisfies $\Phi$ will also satisfy $C$.
\item There exists an $i \notin S_+ \cup S_-$ such that $x_i$ occur in $C$. In that case, there exist $x,x'$ with $x_i\neq x'_i$ that satisfy $\Phi$  (otherwise $x_i$ would have been a spine variable!). Either $x$ or $x'$ will satisfy $C$, so $\Phi\wedge C$ is satisfiable.
\end{enumerate}
So in any case, $\Phi\wedge C$ is satisfiable.
\end{proof}

\subsection*{Overview of proof idea}
In the remainder of this section, we will assume that $p=\frac{1}{2}-b$ for some small positive $b$, and work with $b$ rather than $p$.

Let $P(t,b)$ be the probability that $\Phi_k^{\frac{1}{2}-b}(tn,n)$ (a $(\frac{1}{2}-b)$-biased $k$-CNF) is satisfiable. By studying the partial derivatives of $P$ and estimating the ratio between them, we derive a pair of differential inequalities for the implicit function given by $P(t,b)=\frac{1}{2}$. Solving these inequalities then gives us an upper and a lower bound on the satisfiability threshold.

The $t$-derivative is given by the probability of making the formula unsatisfiable by adding one more clause. \Cref{spinevar:sat} gives us a complete description, in terms of spine variables, of when adding a new clause can turn a satisfiable formula into an unsatisfiable one.

In order to calculate the $b$-derivative, we employ Margulis-Russo's formula from percolation theory (\cref{prop:derivatives}). This will result in something very similar to the $t$-derivative, only depending on the signs of variables slightly differently. To bridge that gap, we study the effects of re-randomizing the signs that spine variables occur with in clauses, conditional on certain other spine variables not being affected (\cref{lemma:spinesupperbound,lemma:spineslowerbound,lemma:spinesconditional}).

Margulis-Russo's formula was proven by Russo~\cite{russosformula} specifically for indicator random variables of increasing events, but the event we are interested in (satisfiability) is not monotone with respect to changing signs. However, in Grimmett's textbook on percolation theory~\cite{grimmettpercolation} there is a generalization of Russo's formula (thm 2.32) to any real-valued random variable.  Here is a version of that theorem, restated with our notation and for finite-dimensional product spaces.
\begin{thm}[Russo's formula, finite case]\label{Russos-formula}
Let $I$ be a finite set, let the probability space $\mathcal{S}=\{-1,1\}^I$ be equipped with the product measure where $\Pr(s_i=-1)=p$ for any $i \in I$, and let $X$ be a real-valued random variable on $\mathcal{S}$. For any $s\in \mathcal{S}$, let $s^{\pm i}$ be $s$ but with the $i$-coordinate set to $\pm 1$. Furthermore, let the \emph{pivotal} $\delta^{i}X$ be defined by
\[
\delta^{i}X(s):=
X(s^{+i})-X(s^{-i}).
\]
Then, for any $0<p<1$, 
\[
\frac{\partial }{\partial p} \E_p[X] = \sum_{i\in I} \E_p[\delta^{i}X].
\]
\end{thm}
\noindent This theorem allows us to calculate the rate of change of $\E_p[X]$ by studying the expected effect of ‘local’ changes to $s$. We will apply it with $X= \mathbf{1}_{\Phi \in \textrm{SAT}}$.

\subsection*{Proof of \cref{thm:parabolicthreshold}}

\begin{prop}
\label{prop:derivatives}
Let $C$ be a $b$-biased random $k$-clause, and let $C_\pm$ be $C$ but with the sign of the first variable changed to $\pm$. Furthermore, let $\rho_\pm := {\Pr\big(\Phi\wedge C_\pm \notin\SAT,\Phi\in\SAT\big)}$ and $\rho :=  \Pr\big(\Phi\wedge C \notin\SAT,\Phi\in\SAT\big)$. Then 
\begin{align*}
\frac{\partial P}{\partial t} &= n\rho = n\left({\Big(\frac{1}{2}+b\Big)\rho_+} + {\Big(\frac{1}{2}-b\Big)\rho_-}\right)
\\
\frac{\partial P}{\partial b} &= tkn(\rho_+-\rho_-)
\end{align*}
\end{prop}
\noindent \begin{proof}
The $t$-derivative is trivial; the prefactor comes from scaling time by $n$. 
For the $b$-derivative, let $X= \mathbf{1}_{\Phi \in \textrm{SAT}}$ and apply Russo's formula. Let $M$ be the set of clauses in $\Phi(t,b)$. 
\begin{equation}\frac{\partial P}{\partial b}=\frac{\partial \E_p[X] }{\partial p} = \sum_{C\in M}\sum_{x\in C} \E_p[\delta^{x,C}X]  , \label{eq:b-derivative}\end{equation}
where $\delta^{x,C}$is the pivot of the variable $x$ in the clause $C$. Now, pick an arbitrary $C\in M$. Let $x$ be the first variable in $C$. The signed pivotal $\delta^{x,C} X$ is $+1$ if $\Phi\wedge C_+$ is unsatisfiable and $\Phi\wedge C_-$ is satisfiable, $-1$ if the reverse holds, and $0$ otherwise. Thus
\begin{align*}
\E_p[\delta^{x,C}X] = 
\quad&\Pr(\Phi\wedge C_+ \notin\SAT,\Phi\wedge C_-\in\SAT )
\\
-&\Pr(\Phi\wedge C_- \notin\SAT,\Phi\wedge C_+\in\SAT )
\\
=\quad&\Pr(\Phi\wedge C_+ \notin\SAT,\Phi\in\SAT )
\\
-&\Pr(\Phi\wedge C_- \notin\SAT,\Phi\in\SAT )
\\
=\quad &\rho_+-\rho_-
\end{align*}
Since every term in the sum in \cref{eq:b-derivative} above have same expected value, we can apply Wald's equation to arrive at 
\[
\sum_{C\in M}\sum_{x\in C}  \E_p[\delta^{x,C}X] =
\E|M| \cdot k\cdot(\rho_+-\rho_-)\]
Noting that $\E |M|=nt$, the proposition follows.
\end{proof}

\begin{lemma}
\label{lemma:derivativeestimates}
For $\rho$, $\rho_+$ and $\rho_-$ defined as in \cref{prop:derivatives} and $\delta$ as in \cref{thm:spinevars},
\begin{align*}
4b\leq \frac{\rho_+-\rho_-}{\rho} \leq \frac{4tk}{\delta}b
\end{align*}
\end{lemma}
\noindent In order to prove \cref{lemma:derivativeestimates}, we will need the estimates in \cref{lemma:spinesupperbound,lemma:spinesconditional,lemma:spineslowerbound}.
\begin{lemma}
\label{lemma:spinesupperbound}
For  $\Phi$ conditioned on being satisfiable and having a spine variable $x$,  the probability that  $x$ is locked to `TRUE' is at most $\frac{1}{2}+\frac{tk}{\delta}b$.
\end{lemma}
\begin{lemma}
\label{lemma:spinesconditional}
Let $\Gamma$ be the event that  $\Phi$ is satisfiable and that the variables $x_1,\ldots x_s$ are spine variables in $\Phi$, of which $x_2,x_3,\ldots x_{s}$ are locked to signs $\sigma_2,\sigma_3, \ldots \sigma_s$ (for some $\sigma_i=\pm 1$). 
Then, conditional on $\Gamma$, the probability that  $x_1$ is locked to `TRUE' is between $\frac{1}{2}-\frac{tk}{\delta}b$ and $\frac{1}{2}+\frac{tk}{\delta}b$.

Furthermore, if we let  $\sigma = \pm 1$ with probability $\frac{1}{2}\pm b$, then (conditional on $\Gamma$) the probability that $x_1$ is locked to $\sigma$ is between $\frac{1}{2}-\frac{tk}{\delta}b^2$ and $\frac{1}{2}+\frac{tk}{\delta}b^2$
\end{lemma}
\begin{lemma}
\label{lemma:spineslowerbound}
For  $\Phi$ conditioned on being satisfiable and having a spine variable $x$,  the probability that  $x$ is locked to `TRUE' is at least $\frac{1}{2}+b$.
\end{lemma}

\noindent \begin{proof}[Proof of \cref{lemma:derivativeestimates}]
Assume (wlog) that the variables in $C$ are $x_1,\ldots x_k$ and they occur with signs $\sigma_1,\ldots \sigma_k$.
Let the events $E_1$,$E_2$, $E_2^-$ and $E_2^+$ be defined in the following way:
\begin{description}
\item[$E_1$:] $\Phi$ is unsatisfiable and the $k$ variables in clause $C_\pm$ are spine variables in $\Phi$.
\item[$E_2$:] These $k$ variables all occur with the opposite sign (in $C$) to the sign they are locked to in $\Phi-C$.
\item[$E^\pm_2$:] These $k$ variables all occur with the opposite sign (in $C_\pm$) to the sign they are locked to in $\Phi-C_\pm$.
\end{description}
First, note that the event $\{\Phi \notin\SAT,\Phi-C\in\SAT\}$ happens if and only if the events $E_1$ and $E_2$ happen, and similarly  for $C_\pm$, $E_1$ and $E_2^\pm$.
Thus $\rho_\pm = \Pr(E_1, E^\pm_2)$ and $\rho = \Pr(E_1, E_2)$. But $\Pr(E_1, E_2) = \Pr(E_2| E_1)\cdot\Pr(E_1)$, and similarly for $\pm$. It follows that  
\begin{align*}
\frac{\rho_+-\rho_-}{\rho} 
=
&\frac{\Pr(E^+_2| E_1)\cdot\Pr(E_1)-\Pr(E^-_2| E_1)\cdot\Pr(E_1)}{\Pr(E_2| E_1)\cdot\Pr(E_1)}
\\
=
&\frac{\Pr(E^+_2| E_1)-\Pr(E^-_2| E_1)}{\Pr(E_2| E_1)}
\end{align*}
We now want to estimate the probabilities $\Pr(E^\pm_2|E_1)$. Starting with $E^+$, we see that
\begin{align*}
\Pr(E^+_2|E_1) = &\Pr(x_1\in S_{-\sigma_1}|E_1)
\\
\cdot &\Pr(x_2\in S_{-\sigma_2}|E_1,x_1\in S_{-\sigma_1})
\\
\cdot &\Pr(x_3\in S_{-\sigma_3}|E_1,x_1\in S_{-\sigma_1},x_2\in S_{-\sigma_2})
\\
&\ldots 
\\
\cdot &\Pr(x_k\in S_{-\sigma_k}|E_1,x_1\in S_{-\sigma_1},\ldots, x_{k-1}\in S_{-\sigma_{k-1}})
\end{align*}
We will use \cref{lemma:spinesupperbound,lemma:spineslowerbound} to estimate the first probability.
\[
\frac{1}{2}+b\leq \Pr(x_1\in S_+|E_1) \leq \frac{1}{2}+\frac{tk}{\delta}b
\] 
For all the subsequent probabilities, we use \cref{lemma:spinesconditional}. 
\[
\frac{1}{2}-\frac{tk}{\delta}b^2
\leq 
\Pr(x_i\in S_{-\sigma_i}|E_1,x_1\in S_{-\sigma_1},\ldots, x_{i-1}\in S_{-\sigma_{i-1}})
\leq
\frac{1}{2}+\frac{tk}{\delta}b^2
\]
For $b$ sufficiently small, the product of these $k-1$ probabilities is at most \[{2^{1-k}(1+ \frac{2tk}{\delta}{b^2})^{k-1}}\leq {2^{1-k}(1+\frac{2tk^2}{\delta}b^2)}\] and similarly at least ${2^{1-k}(1-\frac{2tk^2}{\delta}b^2)}$. Together these estimates yields
\[
2^{-k}\left(1+2b-\frac{2tk^2}{\delta} b^2\right)\leq \Pr(E^+_2|E_1) \leq 2^{-k}\left(1+\frac{2tk^2}{\delta}b+\frac{2tk^2}{\delta} b^2\right)
\]
Similarly, 
\[
2^{-k}\left(1-\frac{2tk^2}{\delta} b-\frac{2tk^2}{\delta} b^2\right)\leq \Pr(E^-_2|E_1) \leq 2^{-k}\left(1-2b+\frac{2tk^2}{\delta} b^2\right)
\]

Using these bounds, we see that 
\begin{align*}
{2^{-k} (4b-O(b^2))}\leq \Pr(E^+_2|E_1)-\Pr(E^-_2|E_1) \leq {2^{-k}( \frac{4tk^2}{\delta}b+O(b^2))}
\end{align*}
Finally, to estimate $\Pr(E_2|E_1)$, simply note that it is a weighted mean of $\Pr(E^+_2|E_1)$ and $\Pr(E^-_2|E_1)$, both of which are $2^{-k}(1\pm O(b))$.
Thus, cancelling the prefactors of $2^{-k}$, we get that 
\[
\frac{ 4b-O(b^2)}{1+ O(b)}
\leq
\frac{\Pr(E^+_2| E_1)-\Pr(E^-_2| E_1)}{\Pr(E_2| E_1)}
\leq
\frac{\frac{4tk^2}{\delta} b-O(b^2)}{1- O(b)}.
\]
The lemma follows.
\end{proof}

\noindent \begin{proof}[Proof of \cref{lemma:spinesupperbound}]
Let $\Gamma$ be the event ${\{\Phi\in \SAT,s(\Phi)\neq 0\}}$. We start by conditioning on $\Gamma$ and on $x$ having degree $d$ in the hypergraph of $\Phi$.
Let $C_1,\ldots, C_d$ be the clauses containing $x$, and let $s_i$ be sign $x$ occurs with in $C_i$. Let $f$ be the function $f:\Sigma_d\to \Sigma_1$ such that $f(\sigma)=1$ iff the formula obtained by replacing $x$ in $C_i$ with $\sigma_i$ for every $i$ is satisfiable (otherwise $f(\sigma)=-1$). Note that $f$ is non-decreasing: satisfying more clauses can only make the rest of the formula easier to satisfy.

When trying to find a satisfying assignment to $\Phi$, we can either satisfy all clauses containing $x$ or all clauses containing $-x$.  
If $f(\mathbf{s})=f(-\mathbf{s})=-1$, then  the formula $\Phi$ is unsatisfiable regardless of the value we assign to $x$, and similarly if $f(\mathbf{s})=f(-\mathbf{s})=1$ it is always satisfiable. But if $f(\mathbf{s})\neq f(-\mathbf{s})$ then $\Phi$ is satisfiable and $x$ is a spine variable locked to $\frac{1}{2}(f(\mathbf{s})-f(-\mathbf{s}))$.

We therefore let $g(\sigma):= \frac{1}{2}(f(\sigma)-f(-\sigma))$. The function $g$ is both non-decreasing and odd.

Now, pick a $b$-biased random $\mathbf{b}\in \Sigma_d$ conditional on $g(\mathbf{b})\neq 0$. Construct a new formula $\Phi'$ from $\Phi$ by replacing the signs $ \mathbf{s}$ that $x$ occur with in $\Phi$ with $\mathbf{b}$. Because of the conditioning, $x$ is a spine variable in $\Phi'$ too, and $\Phi'$ is satisfiable. (Note that $S(\Phi)$ is not necessarily equal to $S(\Phi')$, we only know that $x$ belongs to both.)
What is the expected value of $g(\mathbf{b})$?

Define $W(\sigma)$ to be the probability of $\mathbf{b}=\sigma$, i.e. $W(\sigma) = \left(\frac{1}{2}+b\right)^{h}\left(\frac{1}{2}-b\right)^{d-h}$ where $h$ is the number of $+1$'s in $\sigma$. Furthermore, define $w(g)$ as expected value of $g(\mathbf{b})$, conditional on $g(\mathbf{b})\neq 0$, i.e.
\[
w(g) := \frac{\sum_{\sigma\in \Sigma_d: g(\sigma)=1}  W(\sigma) - W(-\sigma)}{\sum_{\sigma\in \Sigma_d: g(\sigma)=1}   W(\sigma) + W(-\sigma)}.
\]
We will upper bound $\frac{W(\sigma) - W(-\sigma)}{ W(\sigma) + W(-\sigma)}$, and thus get an upper bound for $w(g)$. Cancelling common factors of $W(\sigma) $ and $W(-\sigma)$, we see that
\[
\frac{|W(\sigma) - W(-\sigma)|}{ W(\sigma) + W(-\sigma)}
=
\frac{\left(\frac{1}{2}+b\right)^{|d-2h|}-\left(\frac{1}{2}-b\right)^{|d-2h|}}{\left(\frac{1}{2}+b\right)^{|d-2h|}+\left(\frac{1}{2}-b\right)^{|d-2h|}}.
\]
But for any integer $a$, 
\begin{align*}
\frac{\left(\frac{1}{2}+b\right)^a-\left(\frac{1}{2}-b\right)^a}{\left(\frac{1}{2}+b\right)^a+\left(\frac{1}{2}-b\right)^a}
=
\frac{\sum_{i \textrm{ odd}} \binom{a}{i}(2b)^i}{\sum\limits_{i \textrm{ even}}\binom{a}{i}(2b)^i }
=
2b\cdot \frac{\sum_{i \textrm{ even}} \binom{a}{i-1}(2b)^i}{\sum_{i \textrm{ even}}\binom{a}{i}(2b)^i }.
\end{align*}
Noting that $\binom{a}{i-1}/\binom{a}{i}$ is at most $a$, we see that the above expression is at most $2ab$. Hence the expression $\frac{|W(\sigma) - W(-\sigma)|}{ W(\sigma) + W(-\sigma)}$ is at most ${2b\cdot |d-2h| \leq 2db}$, and it follows that $w(g) \leq 2db$.
So the expected sign of $x$ is at most $2db$, or in other words 
\[
\Pr(x\in S_+|\{\deg(x)=d,  x\in S\}, \Gamma)\leq \frac{1}{2}+db.
\]
We don't know the average degree of spine variables, but we do know the average degree of \emph{all} variables, and that there are at least $\delta n$ spine variables. This gives us the upper bound 
\[
\E[\deg(x)|\{x\in S\}, \Gamma]\leq \frac{\E[ \deg(x)]}{\Pr(x\in S|\Gamma)}.
\]
But $\E [\deg(x)] = tk$, and $\Pr(x\in S|\Gamma) \geq \delta>0$ by \cref{thm:spinevars}. Thus
\[
\Pr(x\in S_+|\deg(x)=d, x\in S,\Gamma)\leq \frac{1}{2}+\frac{tk}{\delta}b.
\]
\end{proof}

\noindent \begin{proof}[Proof of \cref{lemma:spinesconditional}]
We cannot simply re-randomize the signs that $x_1$ appear with, conditional on it remaining a spine variable, because that could change whether or not $x_2$ (say) is a spine variable.

What we can do, however, is re-randomize the signs of $x_1$ conditional on $S_+$ and $S_-$ being unchanged. Consider $D\in \Sigma_d$ of all sign vectors $\sigma=(\sigma_1,\ldots \sigma_d)$ such that replacing the original signs of $x_1$ with $\sigma$ will not change $S_+$ or $S_-$. This set $D$ is symmetric, $D=-D$, so its elements comes in anti-podal pairs.

Now the proof from \cref{lemma:spinesupperbound} carries through as before, giving the upper bound of the corollary. For the lower bound, replace $b$ with $-b$ throughout. 
\end{proof}

Before we continue with the proof of \cref{lemma:spineslowerbound}, we will need the following definition and theorem concerning the possible sizes of simplicial complexes. The theorem was proven independently in \cite{Kruskal,Katona}.
\begin{define}
For positive integers $N$ and $r$, the $r$-\emph{cascade} of $N$ is defined\,\footnote{Our indices here are slightly non-standard; usually one works with $n_{i+1}:=a_i-i$.} as the unique non-increasing sequence of positive integers $a_i$ such that
\[
N=\binom{a_0}{r}+\binom{a_1-1}{r-1}+\ldots +\binom{a_j-j}{r-j}.
\]
Given such $a_i$, we define $N^{(r)}:= \binom{a_0}{r+1}+\binom{a_1-1}{r}+\ldots +\binom{a_j-j}{r-j+1}$. 
\end{define}

\begin{thm}[Kruskal-Katona]
For an integral vector $f$, there exists a $d$-dimensional simplicial complex $\Delta$ such that $f=f_\Delta$ if and only if $0\leq f_r \leq f_{r-1}^{(r)}$ for every $0\leq r\leq d$.
\end{thm}

\noindent\begin{proof}[Proof of \cref{lemma:spineslowerbound}]
Recall that $w(g)$ is the expected sign of a variable, conditional on $g(\mathbf{b})\neq 0$.

Now, $\Delta=g^{-1}(1)$ is a simplicial complex  (equal to  $f^{-1}(1)$ or a subcomplex of it), and this correspondance is a bijection, so $w$ is determined by $\Delta$. Not only that, $w$ only depends on the number of $1$'s of $g$ at each level of the hypercube $\Sigma_d$, so $w$ is determined solely by the $f$-vector $f_\Delta$ of $\Delta$ (the vector whose $i$:th coordinate is the number of $i$-faces of $\Delta$). So henceforth we consider $w$ to be a function of $f_\Delta$.

We want to minimize $w(f)$ over the set of $f$'s that can be written as $f=f_\Delta$ for some $\Delta$. The Kruskal-Katona theorem gives sufficient and necessary conditions for the existence of a simplicial complex with a given $f$-vector.

\begin{claim}
If $f$ is the $f$-vector of some $d$-dimensional simplicial complex $\Delta$ and $f$ minimizes $w(f)$, then $f_r=f^{(r)}_{r-1}$ for every $r$. In other words, there exists a non-increasing sequence $\mathbf{a}=(a_i)_{i=0}^j$ such that $f_r=\sum_{i=0}^j\binom{a_i-i}{r+1}$ for every $r$.
\end{claim}

\noindent\begin{proof} 
Let $F$ be the set of all integral $(d+1)$-vectors $f$ such that $f=f_\Delta$ for some simplicial complex $\Delta$. We want to find $\min_{f\in F} w(f)$.

Note that for every $r$ such that $r > \frac{d}{2}$, increasing $f_r$ by $1$ will decrease $w(f)$ slightly, but if $r < \frac{d}{2}$ doing so will increase it slightly. (For $r=\frac{d}{2}$, changing $f_r$ has no effect on $w(f)$.) So we want to increase $f_r$ for big $r$'s and decrease $f_r$ for small $r$'s, whenever possible.

We therefore let $r_*:=\lceil\frac{d}{2}\rceil$, and consider an arbitrary fixed integer $\ell$ such that $1\leq \ell\leq \binom{d}{r_*}$. Let $F_\ell$ be the set of all $f\in F$ such that $f_{r_*}=\ell$.

We now aim to find the $f$ that achieves $\min_{f\in F_\ell} w(f)$. 
Let $\mathbf{a}$ be the $r_*$-cascade of $\ell$.

By Kruskal-Katona, for any $r>r_*$ we have that $f_r\leq \binom{a_0}{r}+\ldots +\binom{a_j-j}{r-j}$ and that this bound is tight. So we can assume wlog that this holds with equality for all such $r$. 
Similarly, for any $r<r_*$ we have that wlog $f_r= \binom{a_0}{r}+\ldots +\binom{a_j-j}{r-j}$.

So for any $f$ that minimizes $w(f)$ over $F_\ell$, we have that $w(f)$ is determined by the $r_*$-cascade of $f_{r_*}$. It follows that for any non-increasing sequence of non-negative integers $\mathbf{a}$, we can define a simplicial vector $f(\mathbf{a})$ by letting, for each $r$, 
\[f_r(\mathbf{a}): = \binom{a_0}{r}+\ldots +\binom{a_j-j}{r-j}.\]
This will be the unique minimizer (in $F_\ell$) of $w$.
\end{proof}

Since $f(\mathbf{a})$ satisfies Kruskal-Katona (by design), there exists a simplicial complex $\Delta$ (not necessarily unique) with simplicial vector $f(\mathbf{a})$. The size of $\Delta$ is given by

\[
|\Delta|
=
\sum_{r=0}^d \left(\binom{a_0}{r}+\ldots +\binom{a_j-j}{r-j}\right)
=
2^{a_0}+\ldots +2^{a_j-j}
<
2^{a_0+1},
\]
where the last inequality comes from noting that $2^{a_i-i}\leq 2^{a_0-i}$.
But $|\Delta |$ cannot be larger than $2^{d-1}$, because $\Delta$ and $-\Delta$ are disjoint. (By definition, $\Delta=g^{-1}(1)$, and by symmetry $-\Delta:=g^{-1}(-1)$.) 
It follows that either $a_0=d-1$ and $j=0$, or $a_0\leq d-2$. In the former case, $|\Delta|=2^{d-1}$, and in the latter case $|\Delta|< 2^{d-1}$ regardless of the values of $a_1,a_2\ldots a_j$ and $j$.

\begin{claim}
\label{claim:minimizing-w}
For a non-increasing sequence $\mathbf{a} = (a_0,\ldots a_j)$,
\[w(f(\mathbf{a})) = \sum_{i=0}^j \underbrace{\left(\frac{1}{2}-b\right)^{i}\left(\frac{1}{2}+b\right)^{d-a_i}-\left(\frac{1}{2}+b\right)^{i}\left(\frac{1}{2}-b\right)^{d-a_i}}_{=:S(i,a_i)},\]
and $S$ is a strictly decreasing function of both of its arguments.
\end{claim}
\noindent \begin{proof}
Consider a subcube of $\Sigma_d$ whose lowest corner is at level $i$ and whose highest corner is at level $a_i$. For any $r$, this cube will have $\binom{a_i-i}{r-i}$ vertices at level $r$, exactly  matching the $i$-term in the $r$-cascade of $f_r$.
One such cube is 
\[(\underbrace{-1,-1 \ldots -1}_{i\textrm{ times} },\star\ldots \star, \underbrace{1,1 \ldots 1}_{d-a_i \textrm{ times}}).\]
The probability of a $b$-biased random string following the pattern above is precisely  ${(\frac{1}{2}-b)^i(\frac{1}{2}+b)^{d-a_i}}$, so its $w$-weight will be 
\[
S(i,a_i)={\left(\frac{1}{2}-b\right)^i\left(\frac{1}{2}+b\right)^{d-a_i}-\left(\frac{1}{2}+b\right)^i\left(\frac{1}{2}-b\right)^{d-a_i}}.
\]
The claim follows.
\end{proof}

\begin{claim}
The minimum of $w(f(\mathbf{a}))$ is either achieved by $\mathbf{a}=a$ or $\mathbf{a}=\underbrace{(a,\ldots,a)}_{d+1}$ for some integer $a$.
\end{claim}
\noindent \begin{proof}
We will use that $\sum_{i=0}^{j}2^{a_i-i}\leq 2^{d-1}$. If $a_0=d-1$ we must have $\mathbf{a}$ equal to the $1$-term sequence $(d-1)$, because $2^{a_0}=2^{d-1}$ is already as large as the sum can be. So assume instead that $a_0\leq d-2$. 

First, assume $\mathbf{a}$ is not a constant sequence, seeking a contradiction. Then there exists an $i$ with $a_{i-1}>a_{i}$.
Let $\mathbf{a}'$ be $\mathbf{a}$ but with $a_i$ replaced by $a_i+1$. (I.e.\ $\mathbf{a}':=  a_0, \ldots, a_{i-1},a_{i}+1, a_{i+1}, \ldots, a_j $.) The sequence $\mathbf{a}'$ is still non-increasing, and by \cref{claim:minimizing-w}, $w(f(\mathbf{a}'))<w(f(\mathbf{a}))$.

Next, assume $a = a_1=a_2=\ldots=a_j$ for some $0<j<d$, again seeking a contradiction. Then, if $ S(j,a)>0$, decrease the length of $\mathbf{a}$ by $1$, which decreases $w$ by $S(j,a)$. If instead $S(j,a)\leq 0$, increase the length of $\mathbf{a}$ by $1$, which decreases $w$ by
$-S(j+1,a)>-S(j,a)\geq 0$.
\end{proof}

We are now left with only the following candidates for $\mathbf{a}$: $(a)$, for some $a\leq d-1$, or $(a, \ldots ,a)$ for some $a\leq d-2$. It is easy to check that 
$w(f(a)) < w(f(a-1,a-1,\ldots ,a-1 ))$.
Furthermore, $w(f(a)) = (\frac{1}{2}+b)^{d-a}-(\frac{1}{2}-b)^{d-a}$ is decreasing in $a$, so the minimum is achieved for $a = d-1$.
It follows that
\[\min_f w(f)=  \left(\frac{1}{2}+b\right)^{d-(d-1)}-\left(\frac{1}{2}-b\right)^{d-(d-1)}=2b.\]
\end{proof}

\noindent We now have all the ingredients necessary to prove \cref{thm:parabolicthreshold}.

\noindent \begin{proof}[Proof of \cref{thm:parabolicthreshold}]
Let $t=\psi(b)$ be the implicit function defined by $P(t,b)=\frac{1}{2}$.
The derivative of $\psi$ is given by
\[
\psi'(b) = \frac{\partial P / \partial b}{\partial P / \partial t} .
\]
By \cref{prop:derivatives}, $\partial P / \partial b = nkt(\rho_+-\rho_-)$ and  $\partial P / \partial t = n\rho$. So $\psi' = \frac{(\rho_+-\rho_-)}{\rho} k\psi $. But by \cref{lemma:derivativeestimates}, $4b\leq \frac{\rho_+-\rho_-}{\rho}\leq \frac{4tk^2}{\delta}b$, which leads to the following pair of differential inequalities:
\begin{equation}
4k\psi(b)b\leq \psi'(b)\leq \frac{4k^3}{\delta}\psi^2(b)b. \label{diff-ineqs}
\end{equation}
We will relax these inequalities slightly.
First, we see that $\psi'>0$, so $\psi(b)\geq \psi(0)$. By \cref{delta-asymptotics}, $\delta=\delta(t)$ is a decreasing function, and we may replace the function $\delta$ with the constant $\delta_0:=\delta(\psi(0))\geq \exp(-O(k))$ in \cref{diff-ineqs}.

This leads to the differential inequality $\psi'(b)\leq \frac{4k^3}{\delta_0}\psi^2(b)b$, which has the solution $\psi(b)\leq {\psi(0)}\big(1-\frac{2k^3 \psi(0)}{\delta_0}b^2\big)^{-1}$. Noting that $\psi(0)\leq 2^k\log 2$, we see that $\psi(b)< 2^k$ for small enough $b$.
We may therefore also replace $\psi^2$ with $2^k \psi$ in the right-hand side of \cref{diff-ineqs}, leading to 
\begin{equation*}
4k\psi(b)b\leq \psi'(b)\leq \frac{4k^3 2^k}{\delta_0}\psi(b)b,
\end{equation*}
and this new pair of of differential inequalities has the solution
\[
\psi(0)\cdot \exp(2kb^2)\leq \psi(b)\leq \psi(0)\cdot \exp\left(\frac{2k^3 2^k}{\delta_0}b^2\right).
\]
Using the lower bound $\delta_0 \geq \exp(-O(k))$ again,
the theorem follows. 
\end{proof}

\bibliographystyle{abbrv}
\bibliography{satpapers}
\end{document}